\input amssym.def
\input amssym.tex

\def\item#1{\vskip1.3pt\hang\textindent {\rm #1}}


\newskip\litemindent
\litemindent=0.7cm  
\def\Litem#1#2{\par\noindent\hangindent#1\litemindent
\hbox to #1\litemindent{\hfill\hbox to \litemindent
{\ninerm #2 \hfill}}\ignorespaces}
\def\litem{\Litem1}

\tolerance=300
\pretolerance=200
\hfuzz=1pt
\vfuzz=1pt

\hoffset=0in
\voffset=0.5in

\hsize=5.8 true in 
\vsize=9.2 true in
\parindent=25pt
\mathsurround=1pt
\parskip=1pt plus .25pt minus .25pt
\normallineskiplimit=.99pt

\countdef\revised=100
\mathchardef\emptyset="001F 
\chardef\ss="19
\def\3{\ss}
\def\anf{$\lower1.2ex\hbox{"}$}
\def\frac#1#2{{#1 \over #2}}
\def\>{>\!\!>}
\def\<{<\!\!<}

\def\into{\hookrightarrow}
 
\def\ssssarr{\hbox to 15pt{\rightarrowfill}}
\def\sssarr{\hbox to 20pt{\rightarrowfill}}
\def\ssarr{\hbox to 30pt{\rightarrowfill}}
\def\sarr{\hbox to 40pt{\rightarrowfill}}
\def\arr{\hbox to 60pt{\rightarrowfill}}
\def\larr{\hbox to 60pt{\leftarrowfill}}
\def\Arr{\hbox to 80pt{\rightarrowfill}}
\def\mapdown#1{\Big\downarrow\rlap{$\vcenter{\hbox{$\scriptstyle#1$}}$}}

\def\sssmapright#1{\smash{\mathop{\sssarr}\limits^{#1}}}

\def\Alt{\mathop{\rm Alt}\nolimits}

\def\Aut{\mathop{\rm Aut}\nolimits}

\def\Alt{\mathop{\rm Alt}\nolimits}

\def\ch{\mathop{\rm char}\nolimits}

\def\det{\mathop{\rm det}\nolimits}
\def\diag{\mathop{\rm diag}\nolimits}

\def\Ext{\mathop{\rm Ext}\nolimits}

\def\GL{\mathop{\rm GL}\nolimits}

\def\Hom{\mathop{\rm Hom}\nolimits}%
\def\id{\mathop{\rm id}\nolimits} 
\def\im{\mathop{\rm im}\nolimits}


%

\def\ord{\mathop{\rm ord}\nolimits}

\def\rad{\mathop{\rm rad}\nolimits}



\def\SL{\mathop{\rm SL}\nolimits}

\def\span{\mathop{\rm span}\nolimits}

\def\Sp{\mathop{\rm Sp}\nolimits}

\def\supp{\mathop{\rm supp}\nolimits}




\def\0{{\bf 0}}
\def\1{{\bf 1}}

\def\st{{\frak {st}}}

%

\def\K{{{\Bbb K}{\mskip+2mu}}} 

\def\Z{{\Bbb Z}} 
\def\N{{\Bbb N}}

\def\K{{\Bbb K}}

\def\:{\colon}  
\def\.{{\cdot}}
\def\|{\Vert}
\def\bsk{\bigskip}

\def\giantskip{\vskip2\bigskipamount}
\def\gsk{\giantskip}
\def \la {\langle}

\def\msk{\medskip}
\def \ra {\rangle}

\def\bbr{\bigbreak}
\def\giantbreak{\par \ifdim\lastskip<2\bigskipamount \removelastskip
         \penalty-400 \giantskip\fi}

\def\nin{\noindent}
\def\cen{\centerline}
\def\pagebreak{\vskip 0pt plus 0.0001fil\break}
\def\linebreak{\break}

\def\hat{\widehat}

\def\eps{\varepsilon}
\def\epsilon{\varepsilon}

\def\nin{\noindent}
\def\oline{\overline}

\def\pder#1,#2,#3 { {\partial #1 \over \partial #2}(#3)}
\def\pde#1,#2 { {\partial #1 \over \partial #2}}
\def\phi{\varphi}


\def\subeq{\subseteq}

\def\tilde{\widetilde}

\font\ninerm=cmr9
\font\eightrm=cmr8

\font\eightbf=cmbx8


\font\smc=cmcsc10
\font\bfone=cmbx10 scaled\magstep1 
\font\bftwo=cmbx10 scaled\magstep2 

\def\qed{{\unskip\nobreak\hfil\penalty50\hskip .001pt \hbox{}\nobreak\hfil
          \vrule height 1.2ex width 1.1ex depth -.1ex
           \parfillskip=0pt\finalhyphendemerits=0\medbreak}\rm}

\def\qeddis{\eqno{\vrule height 1.2ex width 1.1ex depth -.1ex} $$
                   \medbreak\rm}

\def\Lemma #1. {\bigbreak\vskip-\parskip\noindent{\bf Lemma #1.}\quad\it}

\def\Sublemma #1. {\bigbreak\vskip-\parskip\noindent{\bf Sublemma #1.}\quad\it}

\def\Proposition #1. {\bigbreak\vskip-\parskip\noindent{\bf Proposition #1.}
\quad\it}

\def\Corollary #1. {\bigbreak\vskip-\parskip\nin{\bf Corollary #1.}
\quad\it}

\def\Theorem #1. {\bigbreak\vskip-\parskip\noindent{\bf Theorem #1.}
\quad\it}

\def\Definition #1. {\rm\bigbreak\vskip-\parskip\noindent
{\bf Definition #1.}
\quad}

\def\Remark #1. {\rm\bigbreak\vskip-\parskip\noindent{\bf Remark #1.}\quad}

\def\Example #1. {\rm\bigbreak\vskip-\parskip\noindent{\bf Example #1.}\quad}
\def\Examples #1. {\rm\bigbreak\vskip-\parskip\noindent{\bf Examples #1.}\quad}

\def\Problems #1. {\bigbreak\vskip-\parskip\noindent{\bf Problems #1.}\quad}
\def\Problem #1. {\bigbreak\vskip-\parskip\noindent{\bf Problem #1.}\quad}
\def\Exercise #1. {\bigbreak\vskip-\parskip\noindent{\bf Exercise #1.}\quad}

\def\Conjecture #1. {\bigbreak\vskip-\parskip\noindent{\bf Conjecture #1.}\quad}

\def\Proof#1.{\rm\par\ifdim\lastskip<\bigskipamount\removelastskip\fi\smallskip
            \noindent {\bf Proof.}\quad}

\def\Axiom #1. {\bigbreak\vskip-\parskip\noindent{\bf Axiom #1.}\quad\it}

\def\Satz #1. {\bigbreak\vskip-\parskip\noindent{\bf Satz #1.}\quad\it}

\def\Korollar #1. {\bbr\vskip-\parskip\nin{\bf Korollar #1.} \quad\it}

\def\Folgerung #1. {\bbr\vskip-\parskip\nin{\bf Folgerung #1.} \quad\it}

\def\Folgerungen #1. {\bbr\vskip-\parskip\nin{\bf Folgerungen #1.} \quad\it}

\def\Bemerkung #1. {\rm\bigbreak\vskip-\parskip\noindent{\bf Bemerkung #1.}
\quad}

\def\Beispiel #1. {\rm\bigbreak\vskip-\parskip\noindent{\bf Beispiel #1.}\quad}
\def\Beispiele #1. {\rm\bigbreak\vskip-\parskip\noindent{\bf Beispiele #1.}\quad}
\def\Aufgabe #1. {\rm\bigbreak\vskip-\parskip\noindent{\bf Aufgabe #1.}\quad}
\def\Aufgaben #1. {\rm\bigbreak\vskip-\parskip\noindent{\bf Aufgabe #1.}\quad}

\def\Beweis#1. {\rm\par\ifdim\lastskip<\bigskipamount\removelastskip\fi
           \smallskip\noindent {\bf Beweis.}\quad}

\nopagenumbers

\def\date{\ifcase\month\or January\or February \or March\or April\or May
\or June\or July\or August\or September\or October\or November
\or December\fi\space\number\day, \number\year}

\def\title{Title ??}
\def\author{Author ??}

\def\thanks#1{\footnote*{\eightrm#1}}

\def\rightheadline{\hfil{\eightrm\title}\hfil\tenbf\folio}
\def\leftheadline{\tenbf\folio\hfil{\eightrm\author}\hfil}
\headline={\vbox{\line{\ifodd\pageno\rightheadline\else\leftheadline\fi}}}

\def\firstheadline{}
\def\firstfootline{\cen{\rm\folio}}

\def\seite #1 {\pageno #1
               \headline={\ifnum\pageno=#1 \firstheadline
               \else\ifodd\pageno\rightheadline\else\leftheadline\fi\fi}
               \footline={\ifnum\pageno=#1 \firstfootline\else{}\fi}}

\newdimen\dimenone
 \def\checkleftspace#1#2#3#4{
 \dimenone=\pagetotal
 \advance\dimenone by -\pageshrink   
 \ifdim\dimenone>\pagegoal          
   \else\dimenone=\pagetotal
        \advance\dimenone by \pagestretch
        \ifdim\dimenone<\pagegoal
          \dimenone=\pagetotal
          \advance\dimenone by#1         
          \setbox0=\vbox{#2\parskip=0pt                
                     \hyphenpenalty=10000
                     \rightskip=0pt plus 5em
                     \noindent#3 \vskip#4}    
        \advance\dimenone by\ht0
        \advance\dimenone by 3\baselineskip   
        \ifdim\dimenone>\pagegoal\vfill\eject\fi
          \else\eject\fi\fi}


\def\subheadline #1{\nin\bigbreak\vskip-\lastskip
      \checkleftspace{0.9cm}{\bf}{#1}{\medskipamount}
          \indent\vskip0.7cm\centerline{\bf #1}\medskip}
\def\subsection{\subheadline} 

\def\lsubheadline #1 #2{\bigbreak\vskip-\lastskip
      \checkleftspace{0.9cm}{\bf}{#1}{\bigskipamount}
         \vbox{\vskip0.7cm}\cen{\bf #1}\msk \cen{\bf #2}\bsk}

\def\sectionheadline #1{\bigbreak\vskip-\lastskip
      \checkleftspace{1.1cm}{\bf}{#1}{\bigskipamount}
         \vbox{\vskip1.1cm}\cen{\bfone #1}\bsk}
\def\section{\sectionheadline} 

\def\lsectionheadline #1 #2{\bigbreak\vskip-\lastskip
      \checkleftspace{1.1cm}{\bf}{#1}{\bigskipamount}
         \vbox{\vskip1.1cm}\cen{\bfone #1}\msk \cen{\bfone #2}\bsk}

\def\lchapterheadline #1 #2{\bigbreak\vskip-\lastskip\indent\vskip3cm
                       \cen{\bftwo #1} \msk \cen{\bftwo #2} \gsk}
\def\llsectionheadline #1 #2 #3{\bigbreak\vskip-\lastskip\indent\vskip1.8cm
\cen{\bfone #1} \msk \cen{\bfone #2} \msk \cen{\bfone #3} \nobreak\bsk\nobreak}


\newtoks\literat
\def\[#1 #2\par{\literat={#2\unskip.}%
\hbox{\vtop{\hsize=.15\hsize\nin [#1]\hfill}
\vtop{\hsize=.82\hsize\nin\the\literat}}\par
\vskip.3\baselineskip}

\def\references{
\sectionheadline{\bf References}
\frenchspacing

\entries\par}

\mathchardef\emptyset="001F 
\def\address{Author: \tt$\backslash$def$\backslash$address$\{$??$\}$}

\def\abstract #1{{\narrower\baselineskip=10pt{\noindent
\eightbf Abstract.\quad \eightrm #1 }
\bigskip}}

\def\firstpage{\nin
{\obeylines \parindent 0pt }
\vskip2cm
\centerline{\bfone\title}
\gsk
\centerline{\bf\author}
\vskip1.5cm \rm}

\def\lastpage{\par\vbox{\vskip1cm\nin
\line{
\vtop{\hsize=.5\hsize{\parindent=0pt\baselineskip=10pt\nin\address}}
\hfill} }}

\def\Box #1 { \msk\par\nin 
\centerline{
\vbox{\offinterlineskip
\hrule
\hbox{\vrule\strut\hskip1ex\hfil{\smc#1}\hfill\hskip1ex}
\hrule}\vrule}\msk }

\def\adots{\mathinner{\mkern1mu\raise1pt\vbox{\kern7pt\hbox{.}}
                        \mkern2mu\raise4pt\hbox{.}
                        \mkern2mu\raise7pt\hbox{.}\mkern1mu}}


\pageno=1
\def\stitle{On the classification of rational quantum tori} 
\def\title{\vbox{\centerline{On the classification of rational quantum tori} 
\centerline{and the structure of their automorphism groups}}} 
\def\author{Karl-Hermann Neeb} 
\def\date{14.2.2007}  
\def\rightheadline{\tenbf\folio\hfil{\tt quantor5.tex}\hfil\eightrm\date}
\def\leftheadline{\tenbf\folio\hfil{\rm\stitle}\hfil\eightrm\date}

\def\ord{\mathop{\rm ord}\nolimits}
\def\mod{\mathop{\rm mod}\nolimits}
\def\st{\scriptstyle}

\firstpage 

\abstract{An ${\st n}$-dimensional quantum torus is a twisted group algebra of the 
group ${\st \Z^n}$. It is called rational if all invertible commutators are roots 
of unity. In the present note we describe a normal form for rational 
${\st n}$-dimensional quantum 
tori over any field. Moreover, we show that for 
${\st n = 2}$ the natural exact sequence 
describing the automorphism group of the quantum torus splits over any 
field. \hfill\break 
Keyword: Quantum torus, normal form, automorphisms of quantum tori \hfill\break 
} 

\subheadline{Introduction} 

Let $\K$ be a field and $\Gamma$ an abelian group. 
A {\it $\Gamma$-quantum torus} 
is a $\Gamma$-graded $\K$-algebra $A = \bigoplus_{\gamma \in \Gamma} A_\gamma$, 
for which all grading spaces are one-dimensional 
and all non-zero elements in these spaces are invertible. 
For any basis $(\delta_\gamma)_{\gamma \in \Gamma}$ of such an algebra 
with $\delta_\gamma \in A_\gamma$, we have 
$\delta_\gamma \delta_{\gamma'} = f(\gamma,\gamma') \delta_{\gamma+\gamma'}$, 
where $f \: \Gamma \times \Gamma \to \K^\times$ is a group cocycle. In this 
sense $\Gamma$-quantum tori are the same as {\it twisted group algebras} 
in the terminology of [OP95]. 
Quantum tori arise very naturally in non-commutative geometry 
as non-commutative algebras which are still very close to commutative ones 
(cf.\ [GVF01]) and they also show up in topology (cf.\ [BL04, Sect.~3]). 

For $\Gamma = \Z^n$, we also speak of 
$n$-dimensional quantum tori, also called {\it skew Laurent polynomial rings} 
if the image of $f$ lies in a cyclic subgroup of $\K^\times$ (cf.\ [dCP93]). 
Important special examples arise for $n = 2$ and 
$f(\gamma,\gamma') = q^{\gamma_1 \gamma_2'}$, which leads to an 
algebra $A_q$ with generators $u_1 = \delta_{(1,0)}$ and 
$u_2 = \delta_{(0,1)}$, satisfying 
$u_1 u_2 = q u_2 u_1,$
and their inverses. 
Finite-dimensional quantum tori and their Jordan analogs 
also play a key role in the structure theory of infinite-dimensional 
Lie algebras because they are the natural coordinate structures 
of extended affine Lie algebras ([BGK96], [AABGP97]). 

The first problem we address in this note is the normal form 
of the finite-dimensional {\it rational quantum tori}, i.e., quantum tori 
with grading group $\Gamma = \Z^n$, for which $f$ takes values in the 
torsion group of $\K^\times$. 
Let $P \subeq \K^\times$ be a subset containing for each 
finite element order arising in the multiplicative group $\K^\times$ 
a single representative. We then show in Section~III 
that any rational $n$-dimensional quantum torus $A$ is 
isomorphic to a tensor product 
$$ A \cong A_{q_1} \otimes \cdots \otimes A_{q_{s-1}} 
\otimes A_{q_s^m} \otimes \K[\Z^{n-2s}], 
\leqno(1) $$
where $q_1,\ldots, q_s \in P$ satisfy 
$1 < \ord(q_s)|\ord(q_{s-1})| \cdots |\ord(q_1)$ 
and $\ord(q_s^m) = \ord(q_s)$ (with $m = 1$ if $2s < n$). 
The existence of such a decomposition is not new. 
Under the assumption that the field $\K$ is algebraically closed 
of characteristic zero, (1) can be found in [ABFP05], and a version for skew Laurent 
polynomial rings is stated in [dCP93] (Remark in 7.2). 
Our main new point in Section III are criteria for two such rational 
quantum tori as in (1) to be isomorphic. 

For any $\Z^n$-quantum torus $A$, its group of automorphisms 
is an abelian extension described by a short exact sequence
$$ \1 \to \Hom(\Z^n,\K^\times) \to \Aut(A) \to \Aut(\Z^n,\lambda) \to \1, \leqno(2) $$
where $\lambda \: \Z^n \times \Z^n \to \K^\times, (\gamma, \gamma') 
\mapsto \delta_\gamma \delta_{\gamma'} \delta_\gamma^{-1} \delta_{\gamma'}^{-1}$ 
is the alternating biadditive 
map determined by the commutator map of the unit group $A^\times$, 
and $\Aut(\Z^n,\lambda) \subeq \GL_n(\Z) \cong \Aut(\Z^n)$ is the subgroup 
preserving~$\lambda$. The second main result of this note is that 
for $n = 2$ the sequence (2) always splits. 
In this case $A \cong A_q$ for 
some $q \in \K^\times$, and $\Aut(\Z^2,\lambda) = \GL_2(\Z)$ if 
$q^2 = 1$ and $\Aut(\Z^2,\lambda) = \SL_2(\Z)$ otherwise. 
The statement of 
this result (in case $q$ is not a root of unity) can also be 
found in [KPS94, Th.~1.5], but without any argument for the splitting 
of the exact sequence (2). 
According to [OP95, p.429], the determination of the automorphism groups of 
general quantum tori seems to be a hopeless problem, but we think that 
our splitting result stimulates some hope that more explicit descriptions might 
be possible if the range of the commutator map is sufficiently well-behaved. 

We thank B.~Allison and A.~Pianzola for stimulating discussions on the 
subject matter of this paper, A.~Pianzola for pointing out the 
reference [OP95], B.~Allison for carefully reading earlier versions of 
the manuscript, K.~Goodearl for pointing out references 
[dCP93], [Pa96] and [BL04] and P.~Gille for pointing out 
that the surjectivity of the map $\Phi$ in Proposition~II.3 
can be derived from [Bro82]. Last, but not least, 
we thank the referee for a very detailed report 
that was extremely helpful in improving the exposition of this paper.

\subheadline{Notation} 

Throughout this paper $\K$ denotes an arbitrary field. 
We write $A^\times$ for the unit group of a unital $\K$-algebra $A$. 

Let $\Gamma$ and $Z$ be abelian groups, both written additively. 
A function 
$f \: \Gamma \times \Gamma \to Z$ is called a {\it $2$-cocycle} if 
$$ f(\gamma,\gamma') + f(\gamma + \gamma', \gamma'') 
=  f(\gamma, \gamma'+ \gamma'') + f(\gamma', \gamma'') $$
holds for $\gamma, \gamma', \gamma'' \in \Gamma$. 
The set of all $2$-cocycles is an additive group $Z^2(\Gamma,Z)$  
with respect to pointwise addition. The functions of the form 
$h(\gamma) - h(\gamma+\gamma') + h(\gamma')$
are called {\it coboundaries}. They form a subgroup $B^2(\Gamma,Z) \subeq Z^2(\Gamma,Z)$, 
and the quotient group 
$H^2(\Gamma,Z) := Z^2(\Gamma,Z)/B^2(\Gamma,Z)$ is called the second 
cohomology group of $\Gamma$ with values in $Z$. It classifies central 
extensions of $\Gamma$ by $Z$ up to equivalence. Here we assign to 
$f \in Z^2(\Gamma,Z)$ the central extension 
$Z \times_f \Gamma$, which is the set $Z \times \Gamma$, endowed with the 
group multiplication 
$$ (z,\gamma)(z',\gamma') = (z + z' + f(\gamma, \gamma'), \gamma + \gamma')\quad 
z,z' \in Z, \gamma, \gamma' \in \Gamma. \leqno(0.1) $$

We also write $\Ext(\Gamma,Z) \cong H^2(\Gamma,Z)$ 
for the group of all central extensions of 
$\Gamma$ by $Z$, and 
$\Ext_{\rm ab}(\Gamma,Z)$ for the subgroup corresponding to the 
abelian extensions of the group 
$\Gamma$ by $Z$, which correspond to symmetric $2$-cocycles. 

We call a biadditive map $\Gamma \times \Gamma \to Z$ 
vanishing on the diagonal {\it alternating} and denote the set of these 
maps by $\Alt^2(\Gamma,Z)$. 
A function $q \: \Gamma \to Z$ is called a {\it quadratic form} if the map 
$$\beta_q \: \Gamma \times \Gamma \to Z, \quad (\gamma,\gamma') \mapsto 
q(\gamma + \gamma') - q(\gamma) - q(\gamma') $$
is biadditive. 
Note that we do not require here that $q(n\gamma) = n^2 q(\gamma)$ holds for 
$n \in \Z$ and $\gamma \in \Gamma$. 

For $n\in \N := \{ 1,2,3,\ldots\}$,  
we write $Z[n] := \{ z \in Z \: nz = 0\}$ for the $n$-torsion 
subgroup of $Z$. We also write $\N_0 := \N \cup \{0\}$.

\sectionheadline{I. The correspondence between quantum tori and central extensions} 

\Definition I.1. Let $\Gamma$ be an abelian group. A unital associative $\K$-algebra 
$A$ is said to be a {\it $\Gamma$-quantum torus} if it is $\Gamma$-graded, 
$A = \bigoplus_{\gamma \in \Gamma} A_\gamma,$
with one-dimensional grading spaces $A_\gamma$, and each 
non-zero element of $A_\gamma$ is invertible.\footnote*{In [OP95], these algebras are 
called {\it twisted group algebras.}}

For $\Gamma \cong \Z^d$ we call a $\Gamma$-quantum torus also a {\it $d$-dimensional
  quantum torus}. 
\qed

\Remark I.2. In each $\Gamma$-quantum torus $A$ the set 
$A^\times_h := \bigcup_{\gamma \in \Gamma} \K^\times \delta_\gamma$
of homogeneous units (called {\it trivial units} in [OP95]) 
is a subgroup containing $\K^\times \1 \cong \K^\times$ in its center. 
We thus obtain a central extension  
$$ \1 \to \K^\times \to A^\times_h \to \Gamma \to \1 $$ 
of abelian groups.

It is instructive to see how this can be made more explicit in terms of cocycles, 
which shows in particular that each central extension of $\Gamma$ by 
$\K^\times$ arises as $A^\times_h$ for some $\Gamma$-quantum torus~$A$.  

Let $A$ be a $\Gamma$-quantum torus and pick non-zero 
elements $\delta_\gamma \in A_\gamma$,
so that $(\delta_\gamma)_{\gamma \in \Gamma}$ is a basis of $A$. 
        Then each $\delta_\gamma$ 
is an invertible element of $A$, so that we get 
$$ \delta_\gamma \delta_{\gamma'} = f(\gamma,\gamma') \delta_{\gamma + \gamma'} \quad 
\hbox{ for } \quad \gamma, \gamma' \in \Gamma, \leqno(1.1) $$
where
$f \in Z^2(\Gamma,\K^\times)$ is a $2$-cocycle for which 
$A_h^\times \cong \K^\times \times_f \Gamma$ (cf.~(0.1)). 

Conversely, starting with a cocycle $f \in Z^2(\Gamma,\K^\times)$, we define a 
multiplication on the vector space $A := \bigoplus_{\gamma \in \Gamma} \K \delta_\gamma$ 
with basis $(\delta_\gamma)_{\gamma \in \Gamma}$ 
by 
$\delta_\gamma \delta_{\gamma'} := f(\gamma,\gamma') \delta_{\gamma + \gamma'}.$
Then the cocycle property implies that we get a unital associative algebra, 
and it is clear from the construction that it is a $\Gamma$-quantum torus. 
\qed

\Definition I.3. There 
are two natural equivalence relations between quantum tori. The finest 
one is the notion of {\it graded equivalence:} Two $\Gamma$-quantum tori 
$A$ and $B$ are called {\it graded equivalent} if there is an algebra isomorphism 
$\phi \: A \to B$ with $\phi(A_\gamma) = B_\gamma$ for all $\gamma \in \Gamma$. 

A slightly weaker notion is 
{\it graded isomorphy:} Two $\Gamma$-quantum tori 
$A$ and $B$ are called {\it graded isomorphic} if there is an isomorphism 
$\phi \: A \to B$ and an automorphism $\phi_\Gamma \in \Aut(\Gamma)$ 
with $\phi(A_\gamma) = B_{\phi_\Gamma(\gamma)}$ for all $\gamma \in \Gamma$. 
\qed

The following theorem reduces the corresponding classification problems to 
purely group theoretic ones. 

\Theorem I.4. The graded equivalence classes of $\Gamma$-quantum 
tori are in one-to-one correspondence with the central extensions of the group 
$\Gamma$ by the multiplicative group $\K^\times$, hence 
parametrized by the cohomology group $H^2(\Gamma,\K^\times)$. 

The graded isomorphy classes of $\Gamma$-quantum 
tori are parametrized by the set 
$$ H^2(\Gamma,\K^\times)/\Aut(\Gamma) $$
of orbits of the group 
$\Aut(\Gamma)$ in the cohomology group $H^2(\Gamma,\K^\times)$, where 
the action is given on the level of cocycles 
by $\psi.f := (\psi^{-1})^*f = f \circ (\psi^{-1} \times \psi^{-1})$. 

\Proof. If $\phi \: A \to B$ is a graded equivalence of $\Gamma$-quantum tori, 
then the restriction to the group $A^\times_h$ of homogeneous units 
leads to the commutative diagram 
$$ \matrix{ 
\K^\times & \to & A^\times_h & \to & \Gamma \cr 
\mapdown{\id_{\K^\times}} &  & \mapdown{\phi} & & \mapdown{\id_\Gamma} \cr 
\K^\times & \to & B^\times_h & \to & \Gamma. \cr } $$
This means that the central extensions $A^\times_h$ and $B^\times_h$ of $\Gamma$ 
by $\K^\times$ are equivalent. If, conversely, these extensions are equivalent, 
then any equivalence $\phi \: A^\times_h \to B^\times_h$ extends linearly 
to a graded equivalence $A \to B$. 
Now the observation from Remark~I.2 implies that 
the graded equivalence classes of $\Gamma$-quantum 
tori are parametrized by the cohomology group $H^2(\Gamma,\K^\times) 
\cong \Ext(\Gamma, \K^\times)$. 

If $\phi \: A \to B$ is a graded isomorphism of $\Gamma$-quantum tori, 
then the diagram 
$$ \matrix{ 
\K^\times & \to & A^\times_h & \to & \Gamma \cr 
\mapdown{\id_{\K^\times}} &  & \mapdown{\phi} & & \mapdown{\phi_\Gamma} \cr 
\K^\times & \to & B^\times_h & \to & \Gamma \cr } $$
commutes, which means that the corresponding central extensions 
$A^\times_h$ and $B^\times_h$ are contained in the same orbit of 
$\Aut(\Gamma)$ on $\Ext(\Gamma,\K^\times) \cong H^2(\Gamma,\K^\times)$ (we leave the 
easy verification to the reader). 
Conversely, any isomorphism $\phi \: A^\times_h \to B^\times_h$ of central 
extensions extends linearly to an isomorphism of algebras $A \to B$. 
\qed

\sectionheadline{II. Central extensions of abelian groups} 

In this section $\Gamma$ and $Z$ are abelian groups, written 
additively. We shall derive some 
general facts on the set of equivalence classes $\Ext(\Gamma,Z) \cong H^2(\Gamma,Z)$ 
of central extensions of $\Gamma$ by $Z$. In Sections III and IV below we shall apply these 
to the special case $Z = \K^\times$ for a field $\K$. 

\Remark II.1. Let 
$Z \into \hat \Gamma \sssmapright{q} \Gamma$ be a central extension of the
abelian group $\Gamma$ by the abelian group $Z$ and 
$$\hat \lambda \: \hat \Gamma \times \hat\Gamma \to Z, \quad 
(x,y) \mapsto [x,y] := xyx^{-1}y^{-1} $$
the commutator map of $\hat\Gamma$. Its values lie in $Z$ because 
$\Gamma$ is abelian. Obviously, $\hat \lambda(x, x) =0$, 
and $\hat\lambda$ is an alternating biadditive map (cf.\ [OP95, p.418]). 
Moreover, the commutator map is constant on the fibers of the map $q$, 
hence factors through a biadditive map 
$\lambda \in \Alt^2(\Gamma,Z)$. 

Next we write $\hat\Gamma$ 
as $Z \times_f \Gamma$ with a $2$-cocycle $f \in
Z^2(\Gamma,Z)$. 
For the map $\sigma \: \Gamma \to \hat\Gamma, \gamma \mapsto (0,\gamma)$ 
we then have 
$\sigma(\gamma)\sigma(\gamma') = \sigma(\gamma+\gamma') f(\gamma,\gamma'),$
which leads to 
$$ \eqalign{ \lambda(\gamma,\gamma') 
&= \hat \lambda(\sigma(\gamma),\sigma(\gamma')) 
= \sigma(\gamma)\sigma(\gamma') \big(\sigma(\gamma')\sigma(\gamma)\big)^{-1} \cr
&= \sigma(\gamma+\gamma') f(\gamma,\gamma') \big(\sigma(\gamma+\gamma') f(\gamma',\gamma)\big)^{-1} 
= f(\gamma,\gamma') f(\gamma',\gamma)^{-1}
= f(\gamma,\gamma') - f(\gamma',\gamma). \cr} $$
Therefore the map 
$\lambda_f \in \Alt^2(\Gamma,Z)$ defined by 
$$ \lambda_f(\gamma,\gamma') := 
f(\gamma,\gamma') - f(\gamma',\gamma) \leqno(2.1) $$
can be identified with the commutator map of $\hat\Gamma$. 

Note that the commutator map $\lambda_f$ only depends on 
the cohomology class $[f] \in H^2(\Gamma,Z)$. We thus obtain a group homomorphism 
$$ \Phi \: H^2(\Gamma,Z) \to \Alt^2(\Gamma, Z),\quad 
[f] \mapsto \lambda_f.
\qeddis 

\Remark II.2. Each biadditive map $f \: \Gamma \times \Gamma \to Z$ is a cocycle, 
but it is not true that 
each cohomology class in $H^2(\Gamma,Z)$ has a biadditive representative. 
A typical example is the class corresponding to the exact sequence 
$\0 \to m\Z \to \Z \to \Z/m\Z \to \0$. 
\qed

\Proposition II.3. For abelian groups $\Gamma$ and $Z$ we have a 
split short exact sequence 
$$ \0 \to \Ext_{\rm ab}(\Gamma,Z)\to 
\Ext(\Gamma,Z) \cong H^2(\Gamma,Z) \sssmapright{\Phi} \Alt^2(\Gamma,Z) 
\to \0,$$
describing the kernel of the map $\Phi$. 

\Proof. For the exactness in $\Ext(\Gamma,Z)$, we only have to 
observe that an extension $\hat\Gamma$ of $\Gamma$ by $Z$ 
is an abelian group if and only if the commutator map of $\hat\Gamma$ is 
trivial (cf.\ Remark~II.1). 

The remaining assertions can be found as Exercise~5 in [Bro82, \S V.6]. 
The main point of the argument is to use the short exact 
Universal Coefficient Sequence 
$$ \0 \to \Ext_{\rm ab}(\Gamma,Z)\to 
H^2(\Gamma,Z) \sssmapright{\Psi}  \Hom(H_2(\Gamma),Z) \to \0, $$
then show that $H_2(\Gamma) \cong \Lambda^2(\Gamma)$, which leads to 
an isomorphism 
$\Hom(H_2(\Gamma),Z) \cong \Alt^2(\Gamma,Z)$
([Bro82, Thm.~6.4]), and then to verify that $\Phi$ corresponds to 
$\Psi$ under this identification. 
\qed

In [Bro82], the proof of the surjectivity of $\Phi$ is based on 
the observation that each abelian group is a direct limit of its 
finitely generated subgroups which in turn are products of cyclic groups. 
Below we give a direct argument for the surjectivity of $\Phi$ 
if $\Gamma$ is a direct sum of cyclic groups (the only case relevant 
in the following). We thus obtain 
an explicit description of $H^2(\Gamma,Z)$. 

For the following proposition we recall that, as a consequence of the 
Well-Ordering Theorem, each set $I$ carries a total order. We also recall 
the notation $Z[n] = \{ z \in Z \: nz = 0\}$. 

\Proposition II.4. Let $\Gamma = \bigoplus_{i \in I} \Gamma_i$ 
be a direct sum of cyclic groups 
$\Gamma_i \cong \Z/m_i\Z$, $m_i \in \N_0$. Further let $\leq$ be a total order 
on $I$. Then 
$$ H^2(\Gamma,Z) 
\cong \Ext_{\rm ab}(\Gamma,Z) \oplus \Alt^2(\Gamma,Z)
\cong \prod_{m_i \not=0} Z/m_i Z \oplus 
\prod_{i < j} Z[\gcd(m_i, m_j)], \leqno(2.2) $$
where we put $\gcd(m,0) := m$ for $m \in \N_0$. 
If, in addition, $\Gamma$ is free, then $\Phi$ is an isomorphism, 
$H^2(\Gamma,Z) \cong Z^{\{(i,j) \in I^2\: i < j\}},$
and each cohomology class has a biadditive representative. 

\Proof. To see that $\Phi$ is surjective, let 
$\eta \in \Alt^2(\Gamma,Z)$. If $\gamma_i$ is a generator of $\Gamma_i$, 
we have $\eta(n\gamma_i,m\gamma_i) = nm \eta(\gamma_i, \gamma_i) = 0$ for 
$n,m \in \Z$, so that 
$\eta$ vanishes on $\Gamma_i \times \Gamma_i$. 
We define a biadditive map $f_\eta \: \Gamma \times \Gamma \to Z$ by 
$$ f_\eta(\gamma_i, \gamma_j) := \cases{ 
\eta(\gamma_i, \gamma_j) & for $i > j$, $\gamma_i \in \Gamma_i, \gamma_j \in \Gamma_j$, \cr 
0 & for $i \leq  j$, $\gamma_i \in \Gamma_i, \gamma_j \in \Gamma_j$. \cr} $$
Then $f_\eta$ is biadditive, hence a $2$-cocycle (Remark~II.2), and 
$\Phi(f_\eta) = \eta$. 

Clearly, the assignment $\eta \mapsto f_\eta$ defines an injective homomorphism 
$\Alt^2(\Gamma,Z) \to H^2(\Gamma,Z)$, splitting $\Phi$. 
We know from Proposition~II.3, that $\ker \Phi = \Ext_{\rm ab}(\Gamma,Z)$. 

We next observe that 
$\Alt(\Gamma,Z) 
\cong \prod_{i < j} \Hom(\Gamma_i \otimes \Gamma_j, Z),$
and $\Gamma_i \otimes \Gamma_j \cong 
\Z/\gcd(m_i,m_j) \Z,$
which leads to 
$\Hom(\Gamma_i \otimes \Gamma_j, Z) \cong Z[\gcd(m_i,m_j)].$
On the other hand, 
$$ \Ext_{\rm ab}(\Gamma,Z) 
\cong \prod_{i \in I} \Ext_{\rm ab}(\Gamma_i,Z) 
\cong \prod_{m_i \not=0} Z/m_i Z $$
(cf.\ [Fu70, \S 52]), which leads to (2.2). 

If, in addition, $\Gamma$ is free, then $m_i = 0$ for each $i \in I$, and the assertion 
follows from $\Ext_{\rm ab}(\Gamma,Z) = \0$. 
\qed

\sectionheadline{III. The Normal form of rational quantum tori}

In this section we write $\Gamma := \Z^n$ for the free abelian group of rank $n$. 
For an abelian group $Z$ we write $\Alt_n(Z)$ for the set of {\it alternating 
$(n \times n)$-matrices with entries in $Z$}, i.e., 
$a_{ii} = 0$ for each $i$ and $a_{ij} = - a_{ji}$ for $i \not=j$. 
This is an abelian group with respect to matrix addition. 

Clearly the map $\Alt^2(\Gamma,Z) \to \Alt_n(Z), f \mapsto 
(f(e_i, e_j))_{i,j=1,\ldots, n}$ 
is an isomorphism of abelian groups, so that $\Alt_n(Z) \cong H^2(\Gamma,Z)$ by 
Proposition~II.4. Writing $\lambda_A \in \Alt^2(\Gamma,Z)$ for the alternating 
form $\lambda_A(\alpha,\beta) := \beta^\top A \alpha$ determined by the alternating 
matrix $A$, we have for $g \in \GL_n(\Z) \cong \Aut(\Gamma)$ the relation 
$\lambda_A(g.\alpha,g.\beta) = \beta g^\top A g \alpha,$
so that the orbits of the natural action of $\Aut(\Gamma) \cong \GL_n(\Z)$ on 
the set of alternating forms correspond to the orbits of the action of $\GL_n(\Z)$ 
on $\Alt_n(Z)$ by 
$$ g.A := gA g^\top, \leqno(3.1) $$
where we multiply matrices in $M_n(\Z)$ with matrices in $M_n(Z)$ in the 
obvious fashion. We conclude that 
$$ H^2(\Gamma,Z)/\Aut(\Gamma) \cong \Alt_n(Z)/\GL_n(\Z), \leqno(3.2) $$
the set of $\GL_n(\Z)$-orbits in $\Alt_n(Z)$. 

If $n = n_1 + \ldots + n_r$ is a partition of $n$ and 
$A_i \in M_{n_i}(Z)$, then we write 
$$ A_1 \oplus A_2 \oplus \ldots \oplus A_r 
:= \diag(A_1, \ldots, A_r), $$
for the block diagonal matrix with entries $A_1, \ldots, A_r$. For 
$h_1,\ldots, h_s \in Z$ we further write 
$$ N(h_1,\ldots, h_s) := \pmatrix{0 & h_1 \cr -h_1 & 0 \cr} 
\oplus \pmatrix{0 & h_2 \cr -h_2 & 0 \cr} 
\oplus \ldots \oplus \pmatrix{0 & h_s \cr -h_s & 0 \cr} \oplus \0_{n-2s} 
\in \Alt_n(Z). $$

In the following we shall assume that $Z$ is a cyclic group, hence of the 
form $\Z/(m)$ for some $m \in \N_0$. If $m = 0$, then $Z = \Z$ is a principal 
ideal domain. This is not the case for $m > 0$, 
but $Z$ still carries a natural 
ring structure given by $\oline x \cdot \oline y := \oline{xy}$ for 
$x,y \in \Z, \oline x := x + m\Z$, turning it into a principal ideal ring. 
We write $Z^\times$ for the set of units in $Z$ and note that 
if $m = p_1^{\ell_1}\cdots p_k^{\ell_k}$ is the prime factorization of $m$, then the set 
$$ P := \{ \oline p_1^{j_1}\cdots \oline p_k^{j_k} \: 0 \leq j_i \leq \ell_i, i =1,\ldots, k\} 
\subeq Z $$
is a multiplicatively closed set of representatives for the 
multiplicative cosets of the unit group $Z^\times$. 
We say that $a$ {\it divides} $b$ in $Z$, written 
$a|b$, if $bZ \subeq aZ$. Since each subgroup of $Z$ is cyclic and 
determined by its order, we have 
$$ a | b \quad \Longleftrightarrow \quad \ord(b)|\ord(a)
\quad \hbox{ and } \quad 
aZ = b Z \quad \Longleftrightarrow \quad b \in a Z^\times.$$ 
If $h_1, h_2 \in P$ are non-zero and $h_1|h_2$, then the explicit description of the 
set $P$ shows that there exists a unique element 
$h \in P$ with $h_2 = h_1 h$. We then write $h_2/h_1 := h$. 

Although $Z$ is not a principal ideal domain for $m > 0$, we define 
for a matrix $A \in M_n(Z)$ the {\it determinantal divisor} 
$d_i(A) \in P$, $i =1,\ldots, n$, as the unique element in $P$ 
generating the additive subgroup of $Z$ generated by all $j$-minors 
of the matrix $A$. As a consequence of the Cauchy--Binet Formula 
([New72, II.12]), 
$$d_i(AB) = d_i(BA) = d_i(A) \quad \hbox{ for} \quad  A \in M_n(Z), 
B \in \GL_n(Z). \leqno(3.3) $$
We thus obtain a set of $n$ $P$-valued invariants for the action of 
$\GL_n(Z)$ on $\Alt_n(Z)$ satisfying for $N := N(h_1,\ldots, h_s)$ with 
$h_i \in P$ and $h_1|h_2| \ldots | h_s$: 
$$ d_1(N) = h_1, \quad d_2(N) = h_1^2, \quad 
d_3(N) = h_1^2 h_2, \quad \ldots,\quad d_{2s}(N) = h_1^2 \cdots h_s^2, \quad 
d_j(N) = 0, \ j> 2s. $$
Unfortunately, these invariants do not separate the orbits for  
a {\sl finite} cyclic group $Z$, but they do for $Z = \Z$ 
(cf.\ Theorem~III.2 below and [New72, Th.~II.9]). 

\Theorem III.1. {\rm(Smith normal form over cyclic rings)} 
We consider the action of the group $\GL_n(Z) \times \GL_n(Z)$ on 
$M_n(Z)$ by $(g,h).A := gAh^{-1}$. 
\litem{(1)} Each $\GL_n(Z)^2$-orbit contains a unique matrix of 
the form 
$$ \diag(h_1,\ldots, h_n) \quad \hbox{ with } \quad 
h_i \in P, \ h_1|h_2| \ldots| h_n. $$
\litem{(2)} Each $\SL_n(Z)^2$-orbit contains a unique matrix of 
the form 
$$ \diag(h_1,\ldots, zh_n) \quad \hbox{ with } \quad 
h_i \in P, \ h_1|h_2| \ldots| h_n, z \in Z^\times. $$
\litem{(3)} For $h_i \in P$ with $h_1|h_2|\ldots|h_n$, we 
consider the multiplicative subgroup 
$$ D_{(h_1,\ldots, h_s)} := 
\{ \det(g) \: g \in \GL_n(Z), gN(h_1,\ldots, h_s) g^\top = N(h_1,\ldots, h_s)\} 
\leq Z^\times. $$
Then 
$D_{(h_1,\ldots, h_s)} = Z^\times$ for $2s < n$, and if $2s = n$, then 
$$ \{ z \in Z^\times \: zh_s = h_s\} \subeq 
D_{(h_1,\ldots, h_s)} \subeq \{ z \in Z^\times \: z^2 h_s = h_s\}.$$ 

\Proof. (1) [Br93, Th.~15.24] 

(2) For $z \in Z$ we write 
$\sigma(z) := \diag(1,\ldots, 1, z)$
and observe that $\sigma \: Z^\times \to \GL_n(Z)$ is an embedding, which 
leads to a semidirect product decomposition 
$\GL_n(Z) = \SL_n(Z) \sigma(Z) \cong \SL_n(Z) \rtimes Z^\times.$

{\bf Existence:} 
For each $A \in M_n(Z)$,  (1) implies the 
existence of $g,h \in \GL_n(Z)$ such that 
$N := gAh^{-1} = \diag(h_1, \ldots, h_n)$ as in (1). 
Writing $g = \sigma(z) g_1$ and $h = \sigma(w) h_1$ with $g_1, h_1 \in \SL_n(Z)$, 
it follows that 
$$ g_1 A h_1^{-1} = \sigma(z)^{-1} N \sigma(w) = \sigma(z^{-1}w) N 
= \diag(d_1,\ldots, d_{n-1}, z^{-1}w d_n). $$

{\bf Uniqueness:} Suppose first that $Z = \Z$ is infinite. 
If $d_n = 0$, then there is nothing to show. If $d_n \not=0$, 
then the fact that the determinant function is constant on the 
orbits of $\SL_n(Z)^2$ implies the assertion. 

We may therefore assume that $Z = \Z/(m)$ for some $m > 0$. 
Writing $m = p_1^{m_1} \cdots p_k^{m_k}$ for its prime factorization, 
we obtain a direct product of rings 
$$ \Z/(m) \cong \Z/(p_1^{m_1}) \times \ldots \times \Z/(p_k^{m_k}). $$
For $Z_i := \Z/(p_i^{m_i})$, we accordingly have 
$\SL_n(Z) \cong \prod_{i = 1}^k \SL_n(Z_i)$ and 
$M_n(Z) \cong \prod_{i = 1}^k M_n(Z_i),$
as direct products of groups, resp., rings.
Therefore it suffices to prove the assertion for the case 
$Z = \Z/(p^m)$, where $m \in \N$ and $p$ is a prime. 
Each element $z \in Z$ can be uniquely written as 
$$ z = a_0 + a_1 p + \ldots + a_{m-1} p^{m-1} 
\quad \hbox{ with } \quad 
0 \leq a_i < p. $$
It is a unit if and only if $a_0 \not=0$, i.e., 
$p\not|z$. If $p^k$, $0 \leq k \leq m-1$,  
is the maximal power of $p$ dividing $z$, then 
$$ z 
= \sum_{j = k}^{m-1} a_j p^j 
= p^k \sum_{j = k}^{m-1} a_j p^{j-k}, $$
where the second factor is a unit. Therefore 
$\{ 1,p,p^2,\ldots, p^{m-1}, p^m = 0\}$
is a system of representatives of the multiplicative cosets of $Z^\times$ in 
$Z$. 

{\bf Step 1:} We have to show that if two matrices $D(z_1)$ and $D(z_2)$ of the form 
$$D(z) := \diag(p^{k_1}, \ldots, z p^{k_n}) = \sigma(z) D(1), \quad 
0 \leq k_1 \leq \ldots \leq k_n < m,  $$
lie in the same orbit of $\SL_n(Z)^2$, then $D(z_1) = D(z_2)$. 
Since the orbit of $D(z) = \sigma(z) D(1)$ under 
$\SL_n(Z)^2$ coincides with the set $\sigma(z) (\SL_n(Z)^2.D(1))$, 
it suffices to consider the case $z_1 = 1$ and $z_2 = z \in Z^\times$. 

{\bf Step 2:} We proceed by induction on the size $n$ of the matrices. 
For $n = 1$ the group $\SL_n(Z)$ is trivial, which immediately 
implies the assertion. 

{\bf Step 3:} We reduce the assertion to the 
special case $k_1 = 0$. So let us assume that the 
assertion is correct if $k_1 = 0$ and assume that 
there are $g,h \in \SL_n(Z)$ with 
$g D(1) h = D(z).$
Writing 
$D(z) = p^{k_1} D'(z)$ with 
$D'(z) =  \diag(1, p^{k_2-k_1},\ldots, z p^{k_n-k_1}),$
this means that 
$p^{k_1}(g D'(1) h - D'(z)) = 0,$ 
i.e., that $p^{m-k_1}$ divides each entry of the matrix 
$g D'(1) h - D'(z)$. Over the quotient ring $Z' := \Z/(p^{m-k_1})$ we
then have 
$g D'(1) h = D'(z)$ with $k_1 = 0$. Since we assume that 
the theorem holds in this situation, 
we derive that 
$$z p^{k_n - k_1} \equiv  p^{k_n - k_1} \mod p^{m- k_1}. $$
This means that 
$p^{m-k_1}|(z-1) p^{k_n - k_1}$, and hence that 
$p^{k_1}(z-1) p^{k_n-k_1} = (z-1) p^{k_n} = 0$
in $Z$. 

{\bf Step 4:} Now we consider the special case $k_1 = 0$. 
Let $n_1$ be maximal with $k_{n_1} = 0$, $n_2 := n - n_1$, and write 
elements of $M_n(Z)$ accordingly as $(2 \times 2)$-block matrices. 
We further put $k := k_{n_1 + 1} > 0$. Suppose that 
$g D(1) h = D(z)$ for 
$$ g = \pmatrix{ a & b \cr c & d \cr}, \quad 
h = \pmatrix{ a' & b' \cr c' & d' \cr} \in \SL_n(Z). $$
We write 
$$ D(z) = \pmatrix{ \1 & 0 \cr 0 & p^k D'(z)\cr}, 
\quad \hbox{ where } \quad 
D'(z) = \diag(1, p^{k_{n_1+2}-k}, \ldots, zp^{k_n-k}). $$
If $n_1 = n$, then $D(1) = \1$ is the identity matrix, 
and $z = \det(D(z)) = \det(gD(1)h) = 1$ 
proves the assertion in this case. We may therefore 
assume that $1 \leq n_1 < n$. We now have 
$$ \eqalign{ \pmatrix{ \1 & 0 \cr 0 & p^k D'(z)\cr} 
&=\pmatrix{ a & b \cr c & d \cr}
\pmatrix{ \1 & 0 \cr 0 & p^k D'(1)\cr} 
\pmatrix{ a' & b' \cr c' & d' \cr}\cr
&=\pmatrix{ a a' + p^k  b D'(1) c'& ab' + p^k b D'(1) d' \cr 
ca' + p^k d D'(1) c' & cb' + p^k d D'(1) d'\cr}. \cr} $$
From $aa' + p^k b D'(1) c' = \1$, it follows that 
$aa'\equiv \1 \mod\ p$, hence that $\det(a)\in Z^\times$, 
which means that $a \in \GL_{n_1}(Z)$. Multiplication of 
$g$ from the right with the matrix 
$$ g' := \pmatrix{ a^{-1} & - a^{-1} b \sigma(\det a)\cr 0 & \sigma(\det a)\cr} 
\in \SL_n(Z), $$
leads to the relations
$$ gg' =\pmatrix{ a & b \cr c & d \cr}
 \pmatrix{ a^{-1} & - a^{-1} b \sigma(\det a)\cr 0 & \sigma(\det a)\cr} 
= \pmatrix{ \1 & 0 \cr * & * \cr} $$
and 
$$ \eqalign{ (g')^{-1} D(1) 
&= \pmatrix{ a & b \cr 0 & \sigma(\det a)^{-1}\cr} D(1) 
= \pmatrix{ a & p^k b D'(1) \cr 0 & p^k \sigma(\det a)^{-1} D'(1) \cr} \cr
&= \pmatrix{ a & p^k b D'(1) \cr 0 & p^k D'(1) \sigma(\det a)^{-1} \cr} 
= D(1) \pmatrix{ a & p^k b D'(1) \cr 0 & \sigma(\det a)^{-1} \cr}. \cr}$$
With $h' :=  \pmatrix{ a & p^k b D'(1) \cr 0 & \sigma(\det a)^{-1} \cr},$ 
we thus arrive at 
$$ D(z) = gD(1) h 
= gg' (g')^{-1} D(1) h
= (gg') D(1) (h' h). $$
We may now replace $g$ by $gg'$ and $h$ by $h'h$, 
so that we may assume that $a = \1$ and $b = 0$. 
Now  
$$ \pmatrix{ \1 & 0 \cr 0 & p^k D'(z)\cr} 
=\pmatrix{ \1 & 0 \cr c & d \cr}
\pmatrix{ \1 & 0 \cr 0 & p^k D'(1)\cr} 
\pmatrix{ a' & b' \cr c' & d' \cr}
=\pmatrix{ a' & b'  \cr 
ca' + p^k dc' & cb' + p^k d D'(1) d'\cr} $$
leads to $a' = \1$ and $b' = 0$, which in turn implies 
$p^k dD'(1) d' = p^k D'(z).$
In view of $\det(g) = \det(d) =1$ and 
$\det(h) = \det(d') = 1$, we may now use our induction hypothesis 
that the theorem holds for matrices of smaller size. 
Since we have 
$dD'(1) d' = D'(z)$ in the ring $\Z/(p^{m-k})$, 
we thus obtain $p^{k_n - k}(1 - z) = 0$ modulo $p^{m-k}$. 
This leads to 
$0 = p^k p^{k_n-k}(1-z) = p^{k_n}(1-z)$ modulo $p^m$, and from that 
we derive $D(z) = D(1)$.

(3) If $2s < n$ and $N(z) := N(h_1,\ldots, zh_s)$, then 
$\sigma(z).N(1) = N(z) = N(1)$ for each $z\in Z^\times$, so that 
$D_{(h_1,\ldots, h_s)} = Z^\times$. 

Assume $2s= n$. If $zh_s= h_s$, then $\sigma(z).N(1) = N(1)$ 
implies that $z = \det(\sigma(z)) \in D_{(h_1,\ldots, h_s)}$. 
If, conversely, $z \in D_{(h_1,\ldots, h_s)}$, then we pick 
$g \in \GL_n(Z)$ with $\det(g) = z$ and 
$g.N(1) = N(1)$. Then $\sigma(z)^{-1}g \in \SL_n(Z)$ implies that 
$\sigma(z)^{-1}g.N(1)= N(z^{-1})$ 
lies in the $\SL_n(Z)^2$-orbit of 
$\diag(h_1, h_1,\ldots, h_s, z^{-2} h_s),$
and the assertion follows from (2). For this last argument we use that 
for $w \in Z^\times$, 
$$ \pmatrix{0 & - w^{-1} \cr w & 0 \cr} \in \SL_2(Z) 
\quad \hbox{ satisfies}\quad 
\pmatrix{0 & - w^{-1} \cr w & 0 \cr} 
\pmatrix{0 & w h_s \cr -wh_s & 0 \cr} 
= \pmatrix{h_s & 0\cr 0 & w^2h_s \cr}. 
\qeddis 

\Conjecture III. We believe that if $2s = n$, then  
$zh_s = h_s$ for each $z \in D_{(h_1,\ldots, h_s)}$, so that we have equality in 
Theorem III.1(3), whose present version only implies that 
$D_{(h_1,\ldots, h_s)} \cdot h_s$ can be identified with an elementary 
abelian $2$-group, hence is of cardinality $2^k$ for some $k$. 

The conjecture is true if all $h_i$ coincide. In fact,  
for $h_1 = \ldots = h_s$ and $N := N(h_1,\ldots, h_s)$ we write 
$N = h_s N'$, so that the relation 
$g^\top N g = N$ implies that $h_s\cdot(g^\top N' g - N') = 0$. 
We conclude that $g^\top N' g \equiv N' \mod\ord(h_s)$, 
so that [New72, Th.~VII.21] implies the 
existence of some $\tilde g \in \Sp_{2n}(\Z)$ with 
$g \equiv \tilde g \mod \ord(h_s)$. 
Therefore $\det \tilde g = 1$ implies $\det g \equiv 1 \mod \ord(h_s)$, 
i.e., $\det(g) \cdot h_s = h_s$. 
\qed

The following theorem provides a normal form for the orbits of 
$\GL_n(\Z)$ in $\Alt_n(Z)$ for any cyclic group~$Z$. 
For $Z = \Z$ it follows from Theorem~2.19 in [Pa96]. 

\Theorem III.2. For any cyclic group $Z$ the following assertions hold: 
\litem{(1)} Each $\GL_n(\Z)$-orbit in $\Alt_n(Z)$ contains a matrix of the 
form 
$$ N(h_1,\ldots, h_s),\ 2s < n \quad \hbox{ or } \quad 
N(h_1,\ldots, zh_s), 2s =n, $$
with $z \in Z^\times$ and $0 \not= h_i \in P$ satisfying 
$h_1|h_2| \cdots |h_s$. 
\litem{(2)} If the matrices 
$N(h_1,\ldots, zh_s)$ and $N(h_1',\ldots, z'h_s')$ lie in the same 
$\GL_n(\Z)$-orbit, then $s = s'$ and $h_i' = h_i$ for each $i$. 
\litem{(3)} If $2s < n$ or $Z \cong \Z$, then any corresponding 
$\GL_n(\Z)$-orbit contains 
a unique matrix of the form $N(h_1,\ldots, h_s)$. 
If $2s = n$, then two matrices $N(h_1,\ldots, h_{s-1},zh_s)$ 
and \break $N(h_1,\ldots,  h_{s-1}, wh_s)$ lie in the same orbit if and only if 
$d := zw^{-1} \in \pm D_{(h_1,\ldots, h_s)}$. 
In this case, $d^2 h_s = h_s$  

\Proof. (1) Let $q \: \Z \to Z$ be a surjective homomorphism and 
$q_n \: M_n(\Z) \to M_n(Z)$ the induced homomorphism 
which is equivariant with respect to the action (3.1) of $\GL_n(\Z)$ on both groups. 
If $A \in \Alt_n(Z)$, then its diagonal vanishes and 
$a_{ij} = - a_{ji}$, and there exists a matrix $\tilde A \in \Alt_n(\Z)$ with 
$q_n(\tilde A) = A$. 
As $\Z$ is a principal ideal domain, the Theorem on the Skew Normal Form 
([New72, Thms.~IV.1,IV.2]) implies the existence of 
$g \in \GL_n(\Z)$ with 
$$ g.\tilde A =  N(\tilde h_1, \ldots, \tilde h_t) 
\quad \hbox{ and } \quad \tilde h_1| \tilde h_2 |\cdots |\tilde h_t.$$ 
We then have 
$g.A = q_n(g.\tilde A)
=  N(z_1 h_1, \ldots, z_s h_s),$
where $q(\tilde h_j) = z_j h_j$ with $z_j \in Z^\times$, 
$h_j \in P$, and $s$ is maximal with $h_s \not=0$. 
We further get $h_1|h_2|\cdots|h_s$. 

Next we recall from [New72, Th.~VII.6] that 
$q_n(\SL_n(\Z)) = \SL_n(Z)$, which implies that 
$$ q_n(\GL_n(\Z)) = \{ g \in \GL_n(Z) \: \det g \in \{\pm 1\}\}. \leqno(3.4) $$
For $2s = n$ the matrix 
$$ d := \diag(z_1^{-1},1,z_2^{-1},1,\ldots, z_{s-1}^{-1},1,\ldots, 1, 
z_1\cdots z_{s-1}) \in \SL_n(Z) $$
now satisfies 
$$ d.N(z_1h_1,\ldots, z_s h_s) = N(h_1,\ldots, h_{s-1}, z_1\cdots z_s h_s), $$
and for $2s < n$, the matrix 
$$ d := \diag(z_1^{-1},1,z_2^{-1},1,\ldots, z_{s}^{-1},1,\ldots, 1, 
z_1\cdots z_{s}) \in \SL_n(Z) $$
satisfies 
$$ d.N(z_1h_1,\ldots, z_s h_s) = N(h_1,\ldots, h_s). $$
Since $d \in q(\GL_n(\Z))$, this implies (1). 

(2) The Smith Normal Form of the matrix 
$N(h_1,\ldots, zh_s)$ is 
$\diag(h_1, h_1,\ldots, h_s, h_s, 0,\ldots, 0)$
and for the matrix 
$N(h_1',\ldots, z'h_s')$ we have the normal form 
$\diag(h_1', h_1',\ldots, h_{s'}', h_{s'}', 0,\ldots, 0).$
Therefore Theorem~III.1 implies (2). 

(3) In view of (2), the number $s$ and $h_1,\ldots, h_s$ are uniquely 
determined by the $\GL_n(\Z)$-orbit. 
If $2s < n$, then it follows already that the 
corresponding orbit contains a uniqe matrix of the form $N(h_1,\ldots, h_s)$. 
If $Z = \Z$, then the uniqueness assertion 
follows from the uniqueness of the Skew Normal Form 
([New72, Thms.~IV.1,IV.2]), which follows from 
the fact that the determinantal divisors of 
$N = N(h_1,\ldots, h_s)$ satisfy 
$$ h_1 = d_1(N) = d_2(N)/d_1(N),
\ldots, h_s = d_{2s-1}(N)/d_{2s-2}(N) = d_{2s}(N)/d_{2s-1}(N)$$
and $d_j(N) = 0$ for $j > 2s$. 

It remains to consider the case $2s = n$. 
For $\sigma(z) := \diag(1,\ldots, 1,z)$
and $N(z) := N(h_1,\ldots, zh_s)$, 
we get 
$N(z) = \sigma(z).N(1)$, 
and if there exists a $g \in \GL_n(\Z)$ with 
$g.(\sigma(z).N(1)) = \sigma(w).N(1),$
then 
$\det(\sigma(w)^{-1} g \sigma(z)) = w^{-1}z \det(g) \in D_{(h_1,\ldots, h_s)}.$
In view of $\det(g)\in \{\pm 1\}$, this implies that 
$w^{-1}z \in \pm D_{(h_1,\ldots, h_s)}$. 

If, conversely, 
$w^{-1}z \in \pm D_{(h_1,\ldots, h_s)}$, 
then there exists a matrix $g \in \GL_n(Z)$ fixing 
$N(1)$ with $\det(g) \in \{\pm zw^{-1}\}$. 
Hence 
$\det(\sigma(w)g\sigma(z)^{-1}) \in \{\pm 1\},$ 
and (3.4) imply the existence 
of $g_1 \in \GL_n(\Z)$ with $q_n(g_1) = \sigma(w)g\sigma(z)^{-1}$. We now have 
$$ \eqalign{g_1.N(z) 
&=  g_1\sigma(z).N(1) 
=  \sigma(w)g\sigma(z)^{-1}\sigma(z).N(1)
=  \sigma(w)g.N(1) 
=  \sigma(w)N(1) =  N(w). \cr} 
\qeddis 

\Definition III.3. (a) We call a $\Gamma$-quantum torus {\it rational} if 
the commutator group $C_A$ of $A^\times = A_h^\times$ (cf.\ Proposition~A.1) 
consists of roots of unity in $\K$. We call it {\it of cyclic type} if 
$C_A$ is a cyclic  subgroup of $\K^\times$.  

(b) For each $q \in \K^\times$ we write $A_q$ for the $\Z^2$-quantum torus 
corresponding to the biadditive cocycle $f \: \Z^2 \times \Z^2 \to \K^\times$ determined by 
$$ f(e_1, e_1) = f(e_2, e_2) = f(e_2, e_1) = 1 \quad \hbox{ and } \quad f(e_1, e_2) = q. $$
Then the algebra $A_q$ is generated by $u_1 = \delta_{e_1}$, $u_2 = \delta_{e_2}$  
satisfying 
$u_1 u_2 = q u_2 u_1,$
and their inverses. Then $C_{A_q} = \la q\ra$, so that the 
quantum torus $A_q$ is rational if and only if $q$ is a root of unity.
\qed

\Theorem III.4. {\rm(Normal form of rational quantum tori)} 
Let $\K$ be any field. 

\nin {\rm(a)} For 
any rational $n$-dimensional quantum torus $A$ over $\K$,  the 
commutator group $C_A \subeq \K^\times$ is cyclic. Let 
$q$ be a generator of $C_A$ and choose $P \subeq \Z/(m)$ 
for $m = |C_A| = \ord(q)$ as above. 
Then there exists an $s \in \N_0$ with $2s \leq n$ and 
$h_2|\ldots|h_s$ in $P \setminus \{0\}$  
such that 
$$ A \cong A_{q} \otimes A_{q^{h_2}} \otimes \ldots \otimes A_{q^{h_s}} 
\otimes \K[\Z^{n-2s}] \quad \hbox{ and } \quad 2s < n\leqno(3.5) $$
or 
$$ A \cong A_{q} \otimes A_{q^{h_2}} \otimes \ldots \otimes A_{q^{h_{s-1}}} 
\otimes A_{q^{z h_s}} 
\quad \hbox{ and } \quad 2s = n \leqno(3.6) $$
for some $z \in \N$ with $\ord(q^{zh_s}) = \ord(q^{h_s})$. 

\nin {\rm(b)} If two $n$-dimensional rational quantum tori $A$ and $A'$ 
are (graded) isomorphic, then $C_A = C_{A'}$, both can be described 
by some data $(h_2,\ldots, z h_s)$ and 
$(h_2',\ldots, z' h'_{s'})$ as in {\rm(a)} related to the same choice 
of generator $q$ of $C_A = C_{A'}$. 

\nin{\rm(c)} Two 
$n$-dimensional rational quantum tori $A$ and $A'$ given by such data 
are (graded) isomorphic if and only if 
$s = s'$, $h_i = h_i'$ for $i=2,\ldots, s$, and 
$$z' \in \pm z \cdot D_{(1,h_2,\ldots, h_s)}, $$
where 
$$ D_{(h_1,\ldots, h_s)} = 
\{ \det(g) \: g \in \GL_n(Z), gN(h_1,\ldots, h_s) g^\top = N(h_1,\ldots, h_s)\} 
\leq Z^\times $$
for the ring $Z := \Z/\ord(q)$. 
In this case $z^2 h_s = (z')^2 h_s$ holds in $Z$.  

\Proof. (a) We know from Theorem~I.4 and (3.2) 
that the $\Gamma$-quantum tori over $\K$ 
are classified by the orbits of $\Aut(\Gamma) \cong \GL_n(\Z)$ in 
$H^2(\Gamma,\K^\times) \cong \Alt^2(\Gamma,\K^\times)$. 
In this picture, the rational quantum tori correspond to alternating 
forms $f \in \Alt^2(\Gamma,\K^\times)$ 
on $\Gamma$ whose values are roots of unity. 
Since the group $C_A$ generated by the image of $f$ is generated by 
the finite set 
$f(e_i, e_j)$, $i,j=1,\ldots, n$, it is a finite subgroup of 
$\K^\times$, hence cyclic (cf.~[La93, Th.~IV.1.9]). 

Therefore Theorem~III.2 applies, 
and we see that for $s < 2n$ the quantum torus 
$A$ is isomorphic to one defined by a biadditive 
cocycle $f \: \Gamma \times \Gamma \to C_A \subeq \K^\times$, satisfying 
$(\lambda_f(e_i, e_j))_{i,j} = N(q,q^{h_2}, \ldots, q^{h_s})$, 
where $h_2|\ldots|h_s$. Here $h_1 = 1$ follows from the fact that the 
commutator subgroup of $A^\times$ is generated by $q$. 
The quantum torus $A_f \cong A$ defined by $f$ then satisfies (3.5). 
In the other case we have $2s = n$ and (3.6) holds. 

(b) That (graded) isomorphic quantum tori have the same commutator group is clear. 
Therefore (b) follows from (a).

(c) The remaining assertion now follows from 
Theorem~I.4, combined with Theorem~III.2. 
\qed

\Remark III.5. If $A$ is a $\Z^n$-quantum torus of cyclic type 
and the group of commutators in $A^\times$ is generated by $q \in \K^\times$, 
then the Skew Normal Form over $\Z$ and the argument from the proof of 
Theorem~III.1 imply the existence of 
$h_2|\ldots |h_s \in \N$ and $s \in \N_0$ such that 
$$ A \cong A_{q} \otimes A_{q^{h_2}} \otimes \ldots \otimes A_{q^{h_s}} 
\otimes \K[\Z^{n-2s}] $$
(see the Remark in 7.2 of [dCP93]). 
If $\ord(q) = \infty$, two such decompositions describe isomorphic algebras 
if and only if $s = s'$ and $h_i = h_i'$ for all $i$ (Theorem 2.19 in [Pa96] 
or Theorem~III.2). 
The main point of the preceding theorem is that it gives more precise 
information on the isomorphism classes in the rational case. 
\qed

\sectionheadline{IV. Graded automorphisms of quantum tori} 

In this section we briefly discuss the group of automorphisms of a general 
quantum torus, but our main result only concerns the $2$-dimensional case:
For $A = A_q$ and the corresponding alternating form
$\lambda$ on $\Z^2$, the group $\Aut(A)$ it is a semi-direct product 
$\Hom(\Z^2,\K^\times) \rtimes \Aut(\Z^2,\lambda).$

\Definition IV.1. Let $A$ be a $\Gamma$-quantum torus. We write 
$\Aut_{\rm gr}(A)$ for the group of {\it graded automorphisms of $A$}, 
i.e., all those automorphisms $\phi \in \Aut(A)$ for which there exists an 
automorphism $\phi_\Gamma \in \Aut(\Gamma)$ with 
$\phi(A_\gamma) = A_{\phi_\Gamma(\gamma)}$ for all $\gamma \in \Gamma$.
\qed

Note that Proposition~A.1 in the appendix implies that if $\Gamma$ is torsion free, 
then all units are homogeneous, which implies that each automorphism of $A$ is graded.

\Remark IV.2. We fix a basis $(\delta_\gamma)_{\gamma \in \Gamma}$ of $A$ 
and suppose that $f \in Z^2(\Gamma,Z)$ is the corresponding cocycle determined by (1.1). 
Then for each graded automorphism $\phi$ of $A$ there is an automorphism 
$\phi_\Gamma \in \Aut(\Gamma)$ and a function $\chi \: \Gamma \to \K^\times$ such 
$$ \phi(\delta_\gamma) = \chi(\gamma) 
\delta_{\phi_\Gamma(\gamma)}, \quad \gamma \in \Gamma. \leqno(4.1) $$
Conversely, for a pair $(\chi,\phi_\Gamma)$ of a function 
$\chi \: \Gamma \to \K^\times$ and an automorphism 
$\phi_\Gamma \in \Aut(\Gamma)$ the prescription 
$\phi(\delta_\gamma) := \chi(\gamma) \delta_{\phi_\Gamma(\gamma)}$ 
defines an automorphism of $A$ if and only if 
$$ {(\phi_\Gamma^*f)(\gamma,\gamma') \over f(\gamma,\gamma')} = {\chi(\gamma + \gamma') \over 
\chi(\gamma)\chi(\gamma')} \quad \hbox{ for all } \quad 
\gamma,\gamma' \in \Gamma. \leqno(4.2) $$

Note that if $f$ is biadditive, then $\phi_\Gamma^*f/f$ is biadditive, so that 
$\chi$ is a corresponding $\K^\times$-valued quadratic form. 
If $f$ and $\phi_\Gamma$ are given, then a $\chi$ satisfying (4.2) exists if 
and only if $[\phi_\Gamma^*f] = [f]$ holds in $H^2(\Gamma,Z)$. 
\qed

\Lemma IV.3. The image of the map 
$Q \: \Aut_{\rm gr}(A) \to \Aut(\Gamma), \phi \mapsto \phi_\Gamma$
is the group 
$$ \Aut(\Gamma)_{[f]} := \{ \psi \in \Aut(\Gamma) \: 
[\psi^*f] =[f]\}, $$
which is contained in 
$\Aut(\Gamma,\lambda_f) := \{ \psi \in \Aut(\Gamma) \: \psi^*\lambda_f = \lambda_f\},$
where $\lambda_f(\gamma,\gamma') = {f(\gamma,\gamma')\over f(\gamma',\gamma)}$. 
If, in addition, $\Gamma$ is free, then 
$\Aut(\Gamma)_{[f]} = \Aut(\Gamma,\lambda_f).$

\Proof. Let $\phi_\Gamma \in \Aut(\Gamma)$. 
In view of Remark~IV.2, the existence of $\phi \in \Aut_{\rm gr}(A)$ 
with $Q(\phi) = \phi_\Gamma$ is equivalent to the existence of 
$\chi$ satisfying (4.2), which is equivalent to $[\phi_\Gamma^*f] =[f]$ in 
$H^2(\Gamma,\K^\times)$. 
Since (4.2) implies that $\phi_\Gamma^*f/f$ is symmetric, we have 
$\phi_\Gamma^*\lambda_f = \lambda_{\phi_\Gamma^*f} = \lambda_f.$

If, in addition, $\Gamma$ is free, then Proposition~II.4 entails that 
$\phi_\Gamma^*\lambda_f = \lambda_f$ is equivalent to 
$[\phi_\Gamma^*f] = [f]$ in $H^2(\Gamma,\K^\times)$ 
(cf.\ [OP95, Lemma~3.3(iii)]). 
\qed

From (4.2) we derive in particular that $(\chi,\1)$ defines an automorphism of 
$A$ if and only if $\chi \in \Hom(\Gamma,\K^\times)$, so that we obtain 
the exact sequence 
$$ \1 \to \Hom(\Gamma,\K^\times) \to \Aut_{\rm gr}(A) \to \Aut(\Gamma)_{[f]} \to \1 
\leqno(4.3) $$
(cf.\ [OP95, Lemma~3.3(iii)]). We call the automorphisms of the form $(\chi,\1)$ {\it scalar}. 

\Remark IV.4. If the map $\Phi$ from Proposition~II.4 is not injective, 
then the groups $\Aut(\Gamma,\lambda_f)$ and $\Aut(\Gamma)_{[f]}$ need not coincide, 
but with Proposition~II.3 we obtain a $1$-cocycle 
$$ I \: \Aut(\Gamma,\lambda_f) \to \Ext_{\rm ab}(\Gamma,\K^\times), \quad 
\psi \mapsto [\psi^*f - f] $$
with respect to the right action of $\Aut(\Gamma,\lambda_f)$ on 
$\Ext(\Gamma,\K^\times) \cong H^2(\Gamma,\K^\times)$ by 
$\psi.[f] := [\psi^*f]$. We then have 
$\Aut(\Gamma)_{[f]} = I^{-1}(0).$ 
\qed

In the remainder of this section we restrict our attention to the case, 
where $\Gamma = \Z^n$ is a free abelian group of rank~$n$, which 
implies that $\Aut(\Gamma)_{[f]} = \Aut(\Gamma,\lambda_f)$ and 
that $\Aut(A) = \Aut_{\rm gr}(A)$ (Corollary~A.2). 

\Remark IV.5. (a) For $n = 1$, each alternating biadditive map $\lambda$ on 
$\Gamma$ vanishes, so that $\Aut(\Gamma,\lambda) = \Aut(\Gamma) \cong \{\pm \id_\Gamma\}$. 

(b) For each alternating form $\lambda \: \Gamma \times \Gamma \to \K^\times$ we 
have $-\id_\Gamma \in \Aut(\Gamma,\lambda)$. 

(c) In [OP95], it is shown that if $n \geq 3$ and the subgroup $\la \im(\lambda)\ra$ 
of $\K^\times$ generated by the 
image of $\lambda$ is free of rank ${n \choose 2}$, then 
$\Aut(\Gamma,\lambda_f) = \{\pm \id_\Gamma\}$. 

Moreover, for $n = 3$ and $\la \im(\lambda) \ra$ free of rank $2$,  
[OP95, Prop.~3.7] implies the existence of 
a basis $\gamma_1, \gamma_2, \gamma_3 \in \Gamma$ 
with $\lambda(\gamma_1, \gamma_2) = 1$ and 
$$ \eqalign{ \Aut(\Gamma, \lambda) 
&\cong 
\{ \sigma \in \Aut(\Gamma) \: (\exists a,b \in \Z, \eps \in \{\pm 1\})\ 
\sigma(\gamma_1) = \gamma_1^\eps, \sigma(\gamma_2) = \gamma_2^\eps, \sigma(\gamma_3) 
= \gamma_1^a \gamma_2^b \gamma_3^\eps\}\cr
& \cong \Z^2 \rtimes \{\pm\id_{\Z^2}\}.\cr} $$ 
\qed

We now take a closer look at the case $n = 2$. 
Any alternating form $\lambda \in \Alt^2(\Z^2,\K^\times)$ 
is uniquely determined by $q := \lambda(e_1, e_2)$, which implies 
$\lambda(\gamma,\gamma') = q^{\gamma_1 \gamma_2' - \gamma_2 \gamma_1'}.$
We may therefore assume that a corresponding bimultiplicative 
cocycle $f$ satisfies 
$f(\gamma,\gamma') = q^{\gamma_1 \gamma_2'}$, which leads to the quantum torus 
$A_q$ with two generators $u_i = \delta_{e_i}$ and their inverses, satisfying 
$u_1 u_2 = q u_2 u_1$, as defined in the introduction. 

We start with two simple observations: 

\Lemma IV.6. 
$\Aut(\Z^2, \lambda) 
= \cases{ 
\SL_2(\Z) & for $q^2 \not= 1$ \cr 
\GL_2(\Z) & for $q^2 = 1$. \cr }$ 

\Proof. Clearly $\SL_2(\Z) \subeq \Aut(\Z^2,\lambda) \subeq \GL_2(\Z)$. 
The map $g_0(\gamma) = (\gamma_2,  \gamma_1)$ satisfies 
$\GL_2(\Z) \cong \SL_2(\Z) \rtimes \la g_0 \ra$, and we have 
$$ {g_0^*\lambda(e_1, e_2)\over \lambda(e_1, e_2)} 
=  {\lambda(e_2, e_1)\over \lambda(e_1, e_2)} = q^{-2}. 
\qeddis 

\Example IV.7. (a) On $\Z^2$ the map 
$\chi(\gamma) := \gamma_1 \gamma_2$ is a quadratic form with 
$$ \chi(\gamma+ \gamma') - \chi(\gamma) - \chi(\gamma') 
= \gamma_1 \gamma_2' + \gamma_2 \gamma_1'. $$

(b) On $\Z$ the map $\chi(n) := {n \choose 2}$ is a quadratic form with 
$$ \chi(n+n') - \chi(n) - \chi(n') 
= {(n+ n')(n + n'-1) - n(n-1) - n'(n'-1) \over 2} 
= {nn' + n'n\over 2} = nn'. $$
\qed

From $\SL_2(\Z) \subeq \Aut(\Z^2, \lambda)$, it follows in particular that 
each matrix 
$g = \pmatrix{a & b \cr c & d \cr} \in \SL_2(\Z)$ 
can be lifted to an automorphism of $A_q$.
To determine a corresponding quadratic form $\chi \: \Z^2 \to \K^\times$, 
we have to solve the equation (4.2): 
$$ {(g^*f)(\gamma,\gamma')\over f(\gamma,\gamma')} 
= {\chi(\gamma + \gamma')\over \chi(\gamma)\chi(\gamma')}. $$
The form $g^*f/f$ is determined by its values on the pairs 
$(e_1, e_1), (e_1, e_2)$ and $(e_2,e_2)$: 
$$ (g^*f/f)(e_1, e_1) = f(g.e_1, g.e_1) = q^{ac}, \quad 
(g^*f/f)(e_1, e_2) = f(g.e_1, g.e_2)q^{-1} = q^{ad-1} $$
and 
$(g^*f/f)(e_2, e_2) = f(g.e_2, g.e_2) = q^{bd}.$
This means that 
$$ (g^*f/f)(\gamma,\gamma') 
= q^{ac \gamma_1\gamma_1' + (ad-1)(\gamma_1 \gamma_2' + \gamma_1' \gamma_2) 
+ bd \gamma_2 \gamma_2'}. $$

Before we turn to lifting the full groups $\Aut(\Z^2, \lambda)$ 
to an automorphism group of 
$A$, we discuss certain specific elements of finite order separately. 
\Remark IV.8. (a) For the central element $z = -\1 \in \SL_2(\Z)$, any lift 
$\hat z \in \Aut(A_q)$ is of the form 
$$ \hat z.\delta_\gamma = r^{\gamma_1} s^{\gamma_2} \cdot \delta_{-\gamma} 
\quad \hbox{ for some} \quad r,s \in \K^\times, $$
and any such element satisfies 
$\hat z^2.\delta_\gamma 
= r^{\gamma_1} s^{\gamma_2} \cdot \hat z.\delta_{-\gamma}
= r^{\gamma_1-\gamma_1} s^{\gamma_2-\gamma_2} \cdot \delta_{\gamma}
= \delta_\gamma.$
Hence each lift $\hat z$ of $z$ is an element of order $2$. 

(b) The matrices 
$$ g_1 := \pmatrix{ 0 & 1 \cr -1 & 0 \cr} 
\quad \hbox{ and } \quad 
g_2 := \pmatrix{ 1 & 1 \cr -1 & 0 \cr} $$
satisfy $g_1^2 = z= g_2^3$, which leads to 
$\ord(g_1) = 4$ and $\ord(g_2) = 6$. 
From the preceding paragraph we conclude that for any lift 
$\hat g_j$ of $g_j$, $j = 1,2$, we have 
$\hat g_1^4 = \1 = \hat g_2^6.$ 

In view of 
$$ (g_1^*f/f)(\gamma,\gamma') 
= q^{-(\gamma_1 \gamma_2' + \gamma_1' \gamma_2)}, $$
a lift $\tilde g_1$ of $g_1$ is given by 
$\tilde g_1.\delta_\gamma =  q^{-\gamma_1\gamma_2} \delta_{g_1.\gamma}$ 
(Example~IV.7(a)). 
We then have 
$$ \tilde g_1^2.\delta_\gamma 
= q^{-\gamma_1 \gamma_2}\tilde g_1.\delta_{(\gamma_2, -\gamma_1)}
= q^{-\gamma_1 \gamma_2} q^{\gamma_2 \gamma_1}.\delta_{-\gamma}
= \delta_{-\gamma}. $$
Any other lift $\hat g_1$ of $g_1$ is of the form 
$$ \hat g_1.\delta_g = r_1^{\gamma_1} s_1^{\gamma_2} 
q^{-\gamma_1\gamma_2} \delta_{g_1.\gamma} $$
for two elements $r_1, s_1 \in \K^\times$. The square of this element is given by 
$$ \hat g_1^2.\delta_g 
= r_1^{\gamma_1} s_1^{\gamma_2} \hat g_1 \tilde g_1.\delta_{\gamma} 
= r_1^{\gamma_1+\gamma_2} s_1^{\gamma_2-\gamma_1} \tilde g_1^2.\delta_{\gamma} 
= \Big({r_1\over s_1}\Big)^{\gamma_1} 
(r_1 s_1)^{\gamma_2} \cdot \delta_{-\gamma}. \leqno(4.4) $$

For the matrix $g_2$ we have 
$$ (g_2^*f/f)(\gamma,\gamma') 
= q^{-\gamma_1\gamma_1' - (\gamma_1 \gamma_2' + \gamma_1' \gamma_2)}, $$
so that we obtain a lift $\tilde g_2$ of $g_2$ by 
$\tilde g_2.\delta_\gamma = q^{-{\gamma_1 \choose 2} - \gamma_1\gamma_2} 
\delta_{(\gamma_1 + \gamma_2, - \gamma_1)}$
(Example~IV.7(b)). 
Hence each lift $\hat g_2$ of $g_2$ is of the form 
$$ \hat g_2.\delta_g = r_2^{\gamma_1} s_2^{\gamma_2} 
q^{-{\gamma_1 \choose 2} - \gamma_1\gamma_2} 
\delta_{(\gamma_1 + \gamma_2, - \gamma_1)}, $$
for some $r_2, s_2 \in \K^\times$. 
In view of 
$g_2^2 = \pmatrix{0 & 1 \cr - 1 & -1 \cr},$
we get with Example~IV.7(b): 
$$ \eqalign{ \tilde g_2^3.\delta_\gamma 
&= q^{-{\gamma_1 \choose 2} -\gamma_1 \gamma_2} \tilde g_2^2.
\delta_{\gamma_1 +\gamma_2,-\gamma_1} 
= q^{-{\gamma_1 \choose 2} -\gamma_1 \gamma_2} 
q^{-{\gamma_1 + \gamma_2 \choose 2} + (\gamma_1 +\gamma_2)\gamma_1} 
\tilde g_2.\delta_{\gamma_2,-\gamma_1-\gamma_2} \cr 
&= q^{-2{\gamma_1 \choose 2} -{\gamma_2 \choose 2} - \gamma_1 \gamma_2 + \gamma_1^2} 
q^{-{\gamma_2 \choose 2} + (\gamma_1 +\gamma_2)\gamma_2}
\delta_{-\gamma} 
= q^{-\gamma_1(\gamma_1-1) -\gamma_2(\gamma_2-1)  + \gamma_1^2 + \gamma_2^2} 
\delta_{-\gamma} 
= q^{\gamma_1+ \gamma_2} \delta_{-\gamma}. \cr} $$
This further leads to 
$$ \leqalignno{ \hat g_2^3.\delta_\gamma 
&= r_2^{\gamma_1} s_2^{\gamma_2} \hat g_2^2 \tilde g_2.\delta_\gamma 
= r_2^{2\gamma_1 + \gamma_2} s_2^{-\gamma_1 + \gamma_2} \hat g_2 \tilde g_2^2.\delta_\gamma 
= r_2^{2\gamma_1 + 2\gamma_2} s_2^{-2\gamma_1} \tilde g_2^3.\delta_\gamma \cr
&= r_2^{2(\gamma_1 + \gamma_2)} s_2^{-2\gamma_1} q^{\gamma_1 + \gamma_2}.
\delta_{-\gamma}
= \Big({r_2^2\over s_2^2}q\Big)^{\gamma_1} (r_2^2 q)^{\gamma_2} \delta_{-\gamma}. 
&(4.5) \cr}$$

(c) If, in addition, $q^2 = 1$, then 
$\Aut(\Gamma, \lambda_f) = \Aut(\Gamma) \cong \GL_2(\Z)$ 
(Lemma~IV.6). For the involution 
$$g_0 := \pmatrix{0 & 1 \cr 1 & 0\cr}$$
we have $\GL_2(\Z) = \SL_2(\Z) \rtimes \la g_0\ra,$ 
and the elements $g_0, g_1, g_2$ satisfy 
$$ g_0 g_1 g_0 = g_1^{-1} =g_1^3 
\quad \hbox{ and } \quad 
g_0 g_2 g_0 = g_2^5 = g_2^{-1}. \leqno(4.6) $$

To lift $g_0$ to an automorphism of $A_q$, we first note 
that $q^2=1$ implies that 
$$ (g_0^*f/f)(\gamma,\gamma') 
= q^{\gamma_2 \gamma_1' - \gamma_1 \gamma_2'}, 
= q^{\gamma_2 \gamma_1' + \gamma_1 \gamma_2'}, $$
which shows that each lift $\hat g_0$ of $g_0$ is of the form 
$\hat g_0.\delta_\gamma = r_0^{\gamma_1} s_0^{\gamma_2} 
q^{\gamma_1 \gamma_2} \delta_{(\gamma_2, \gamma_1)}$
for some $r_0, s_0 \in \K^\times$. 
In view of 
$$ \hat g_0^2.\delta_\gamma 
= r_0^{\gamma_1} s_0^{\gamma_2} q^{\gamma_1 \gamma_2} \hat g_0.\delta_{(\gamma_2, \gamma_1)} 
= r_0^{\gamma_1+\gamma_2} s_0^{\gamma_2 + \gamma_1} q^{2\gamma_1 \gamma_2} \delta_\gamma 
= (r_0s_0)^{\gamma_1+\gamma_2} \delta_\gamma, $$
$\hat g_0^2 = \1$ is equivalent to 
$r_0 s_0 = 1.$ If this condition is satisfied, then 
$\hat g_0.\delta_\gamma = r_0^{\gamma_1-\gamma_2} 
q^{\gamma_1 \gamma_2} \delta_{(\gamma_2, \gamma_1)}.$ 
\qed

Before we state the following theorem, we recall that for any split extension 
$$ \1 \to A \to \hat G \sssmapright{q} G \to \1 $$
of a group $G$ by some (abelian) $G$-module $A$, the set of 
all splittings is parametrized by the group 
$$Z^1(G,A) = \{ f \: G \to A \: (\forall x,y\in G)\ f(xy) = f(x) +x.f(y)\} $$
of $A$-valued $1$-cocycles. This parametrization is obtained by choosing a 
homomorphic section $\sigma_0 \: G \to \hat G$ and then 
observing that any other homomorphic section 
$\sigma \: G \to \hat G$ is of the form $\sigma = f \cdot \sigma_0$, where 
$f \in Z^1(G,A)$. 

\Theorem IV.9. For each element $q \in \K^\times$ and 
$\lambda(\gamma, \gamma') = q^{\gamma_1 \gamma_2' - \gamma_2 \gamma_1'}$ 
the exact sequence 
$$ \1 \to \Hom(\Z^2, \K^\times) \to \Aut(A_q) \to \Aut(\Z^2, \lambda) \to \1 $$ 
splits. 
For $q^2 = 1$, the homomorphisms 
$\sigma \: \GL_2(\Z) \to \Aut(A_q)$ splitting the 
sequence are parametrized by the abelian group 
$$ Z^1(\GL_2(\Z), \Hom(\Z^2, \K^\times)) \cong 
\{ (r_0, r_1, r_2) \in (\K^\times)^3 \:  r_2^4  r_0^2 = r_1^2\}, $$ 
and for $q^2 \not= 1$, the homomorphisms $\sigma \: \SL_2(\Z) \to \Aut(A_q)$ splitting the 
sequence are parametrized by 
$$ Z^1(\SL_2(\Z), \Hom(\Z^2, \K^\times)) \cong 
(\K^\times)^2 \times \{z \in \K^\times \: z^2 = 1\}. $$

\Proof. First we consider the case $q^2 \not= 1$, 
where $\Aut(\Z^2, \lambda) = \SL_2(\Z)$ 
(Remark~IV.8). We shall use the description of the lifts of 
$g_1, g_2$ given in Remark~IV.8. 
Since $\SL_2(\Z)$ is presented by the relations 
$$ g_1^4 = g_2^6 = \1, \quad g_1^2 = g_2^3 $$
([Ha00, p.51]), Remark~IV.8 implies that 
a pair of elements $(\hat g_1, \hat g_2)$ lifting $(g_1, g_2)$ leads to a lift 
$\SL_2(\Z) \to \Aut(A_q)$ if and only if 
$\hat g_1^2 = \hat g_2^3$. 
Comparing (4.4) and (4.5), we see that $\hat g_1^2 = \hat g_2^3$ is equivalent to 
$$ {r_1\over s_1} = {r_2^2 \over s_2^2} q 
\quad \hbox{ and } \quad 
r_1s_1 = r_2^2 q,  $$
 which is equivalent to 
$$ s_1^2 = s_2^2 
\quad \hbox{ and } \quad 
s_1 = {r_2^2 q \over r_1}, \leqno(4.8) $$
These equations have the simple solution 
$r_1 = q, r_2 = s_1 = s_2 = 1,$
showing that the action of the group $\SL_2(\Z)$ on $\Gamma$ lifts to an action on $A_q$. Moreover, for each pair $(r_1, r_2)$, the set of all 
solutions is determined by the choice of sign in 
$s_2 := \pm s_1$, which is vacuous if $\ch(\K)= 2$. 

\msk

Next we consider the case $q^2 = 1$. 
We assume that the lift $\hat g_0$ of $g_0$ satisfies 
$\hat g_0^2 = \1$ (cf.\ Remark~IV.8(c)). 
Now the 
relation $\hat g_0 \hat g_1 \hat g_0 = \hat g_1^{-1}$ is equivalent to 
$(\hat g_0 \hat g_1)^2 = \1$. We calculate 
$$ \hat g_0 \hat g_1.\delta_\gamma 
= r_1^{\gamma_1} s_1^{\gamma_2} q^{-\gamma_1\gamma_2} \hat g_0.\delta_{(\gamma_2, - \gamma_1)}
= (r_0 r_1)^{\gamma_1} (r_0 s_1)^{\gamma_2} \delta_{(-\gamma_1, \gamma_2)} $$
to get 
$$ (\hat g_0 \hat g_1)^2.\delta_\gamma 
= (r_0 r_1)^{\gamma_1} (r_0 s_1)^{\gamma_2} \hat g_0 \hat g_1.\delta_{(-\gamma_1, \gamma_2)}
= (r_0 s_1)^{2\gamma_2} \delta_\gamma. $$
Hence $\hat g_0 \hat g_1 \hat g_0 = \hat g_1^{-1}$ is equivalent to 
$$ r_0^2 s_1^2 = 1. \leqno(4.9) $$

To see when $\hat g_0 \hat g_2 \hat g_0 = \hat g_2^{-1}$ holds, 
we first observe that 
$$ \hat g_2^{-1}.\delta_\gamma 
= r_2^{\gamma_2} s_2^{-\gamma_1 -\gamma_2} q^{{-\gamma_2 \choose 2}- \gamma_2(\gamma_1 + \gamma_2)} 
\delta_{(-\gamma_2, \gamma_1 + \gamma_2)}. $$
Further 
$$ \eqalign{ \hat g_0 \hat g_2.\delta_\gamma 
&= r_2^{\gamma_1} s_2^{\gamma_2} q^{-{\gamma_1 \choose 2} - \gamma_1\gamma_2} \hat g_0.
\delta_{(\gamma_1 + \gamma_2, - \gamma_1)} 
= (r_0^2 r_2)^{\gamma_1} (r_0 s_2)^{\gamma_2} q^{{\gamma_1 \choose 2} + \gamma_1\gamma_2+ 
(\gamma_1 + \gamma_2)\gamma_1} \cdot \delta_{(-\gamma_1, \gamma_1 + \gamma_2)}\cr 
&= (r_0^2 r_2)^{\gamma_1} (r_0 s_2)^{\gamma_2} q^{{\gamma_1 \choose 2} + 
\gamma_1^2}\cdot\delta_{(-\gamma_1, \gamma_1 + \gamma_2)} 
= (r_0^2 r_2 q)^{\gamma_1} (r_0 s_2)^{\gamma_2} q^{{\gamma_1 \choose 2}}\cdot
\delta_{(-\gamma_1, \gamma_1 + \gamma_2)} \cr} $$
because $q^2 = 1$ implies $q^{n^2} = q^n = q^{-n}$ for each $n \in \Z$. 
 
On the other hand, we have 
$$ \eqalign{ \hat g_2^{-1} \hat g_0.\delta_\gamma 
&= r_0^{\gamma_1} r_0^{-\gamma_2} q^{\gamma_1\gamma_2} \hat g_2^{-1}.\delta_{(\gamma_2, \gamma_1)} 
= r_0^{\gamma_1} r_0^{-\gamma_2} q^{\gamma_1\gamma_2} 
r_2^{\gamma_1} s_2^{-\gamma_2 -\gamma_1} q^{{-\gamma_1 \choose 2}- \gamma_1(\gamma_1 + \gamma_2)} 
\delta_{(-\gamma_1, \gamma_1 + \gamma_2)} \cr 
&= (r_0 r_2 s_2^{-1})^{\gamma_1} (r_0 s_2)^{-\gamma_2}  q^{{-\gamma_1 \choose 2}- \gamma_1^2} 
\delta_{(-\gamma_1, \gamma_1 + \gamma_2)} 
= (r_0 r_2 s_2^{-1})^{\gamma_1} (r_0 s_2)^{-\gamma_2}  q^{-{\gamma_1 \choose 2}} 
\delta_{(-\gamma_1, \gamma_1 + \gamma_2)} \cr 
&= (r_0 r_2 s_2^{-1})^{\gamma_1} (r_0 s_2)^{-\gamma_2}  q^{{\gamma_1 \choose 2}} 
\delta_{(-\gamma_1, \gamma_1 + \gamma_2)}. \cr} $$
Therefore $\hat g_0 \hat g_2 \hat g_0 = \hat g_2^{-1}$ is equivalent to 
$r_0 r_2 s_2^{-1} = r_0^2 r_2 q$ and $(r_0 s_2)^2 = 1,$
which is equivalent to 
$$ r_0 s_2 = q, \leqno(4.10) $$
because this relation implies $(r_0 s_2)^2 = q^2 = 1$. 

We conclude that the numbers $r_0, r_1, r_2, s_1, s_2$ which determine 
$\hat g_0, \hat g_1, \hat g_2$ define a lift of $\GL_2(\Z)$ to 
$\Aut(A_q)$ if and only if the equations (4.8), (4.9) and (4.10) are satisfied: 
$$ s_1^2 = s_2^2, \quad 
s_1 = {r_2^2 q \over r_1}, \quad 
 r_0^2 s_1^2 = 1, \quad \hbox{ and } \quad 
r_0 s_2 = q. $$
If $r_0, r_1$ and $r_2$ are given, we determine 
$s_1$ and $s_2$ by 
$s_1 := {r_2^2 q \over r_1}$ and $s_2 := {q \over r_0}.$
Then 
$$ {s_1^2\over s_2^2} = {r^4_2  r_0^2\over r_1^2} = r_0^2 s_1^2, $$
so that we obtain only the relation 
$r_2^4  r_0^2 = r_1^2$ for $r_0, r_1, r_2$. 
This completes the proof. 
\qed

\Remark IV.10. (a) From the 
proof of the preceding theorem, we see that if $q^2 = 1$, we obtain the particularly simple solution 
$$ r_0 = r_1 = r_2 = 1, \quad s_1 = s_2 = q. $$

(b) For $\ch \K = 2$ the equation $q^2 = 1$ has the unique solution 
$q = 1$, so that $A_q \cong \K[\Z^2]$, and the action of $\GL_2(\Z)$ has a canonical 
lift to an action on $A_q$. 
\qed

\Problem IV.1. Does the sequence (4.3) always split? 
We have seen above, that this is true for 
$\Gamma = \Z^2$. 
If the answer is no, it would be of some interest to understand the 
cohomology groups 
$$H^2(\Aut(\Gamma)_{[f]}, \Hom(\Gamma,\K^\times)) $$
parametrizing the possible abelian extensions of $\Aut(\Gamma)_{[f]}$ by the module 
$\Hom(\Gamma,\K^\times)$. 
\qed

\Problem IV.2. Let $\lambda\in \Alt^2(\Z^n,Z)$, where $Z$ is a cyclic group. 
Determine the structure of the group 
$\Aut(\Z^n,\lambda)$. It should have a semidirect product structure, 
where the normal subgroup is something like a Heisenberg group 
and the quotient is the automorphism group of $\Z^n/\rad(\lambda)$, 
endowed with the induced non-degenerate form. Can this group 
be described in a conventient way by generators and relations? 
Maybe the results in [Is03] can be used to deal with degenerate cocycles. 
\qed

\sectionheadline{A. The group of units if $\Gamma$ is torsion free} 

The following result is used in [OP95, Lemma~3.1] without reference. Here we provide 
a detailed proof. 
\Proposition A.1. If the group $\Gamma$ is torsion free and $A$ a $\Gamma$-quantum torus, 
then $A^\times = A^\times_h$, i.e., each unit of $A$ is graded. 

\Proof. Let $a \in A^\times$ be a unit and 
write $a = \sum_\gamma a_\gamma \delta_\gamma$ in terms of some graded 
basis. We do the same with its inverse 
$a^{-1} = \sum_\gamma (a^{-1})_\gamma \delta_\gamma$, and observe that the set 
$\supp(a) := \{ \gamma \in \Gamma \: a_\gamma \not=0 \}$
is finite. The same holds for $\supp(a^{-1})$, so that both sets 
generate a free subgroup $F$ of $\Gamma$. 
Then $A_F := \span \{ \delta_\gamma \: \gamma \in F\}$ 
is an $F$-quantum torus with $a \in A_F^\times$. We may therefore assume 
that $\Gamma = \Z^d$ for some $d \in \N_0$. 

We prove by induction on $k \in \{0,\ldots, d\}$ that the subalgebra 
$A_k := \span \{ \delta_\gamma \: \gamma \in \Z^{k}\times \{0\}\}$
has no zero-divisors (cf.\ Th.~1.2 in [Pa96]) 
and that all its units are homogeneous. This holds trivially for $k = 0$. 

Let $u_i := \delta_{e_i}$, where $e_1, \ldots, e_d$ is the canonical basis of 
$\Z^d$. We write 
$0 \not= x \in A$ as a finite sum $\sum_{k = k_0}^{k_1} x_k u_d^k$ with 
$x_k \in A_{d-1}$ and $x_{k_0}$ and $x_{k_1}$ non-zero. 
Likewise we write $0 \not= y \in A$ as 
$\sum_{m = m_0}^{m_1} y_m u_d^m$ with 
$y_m \in A_{d-1}$ and $y_{m_0}$ and $y_{m_1}$ non-zero. 
Then the lowest degree term with respect to $u_d$ in $xy$
is 
$$  x_{k_0} u_d^{k_0} y_{m_0} u_d^{m_0} 
= x_{k_0} \big(u_d^{k_0} y_{m_0} u_d^{-k_0}\big) u_d^{k_0+m_0}, $$
and the induction hypothesis implies 
$x_{k_0} u_d^{k_0} y_{m_0} u_d^{-k_0} \not=0$
because conjugation with $u_d$ preserves the subalgebra $A_{d-1}$. 
This implies that $xy\not=0$. 

Now assume that $x \in A$ is a unit and $y = x^{-1}$. 
Since $A_{d-1}$ has no zero-divisors, 
$$ x_{k_0} u_d^{k_0} y_{m_0} u_d^{-k_0} \in A_{d-1}\setminus \{0\} $$
leads to $k_0 + m_0 = 0$. 
A similar consideration for the highest order term implies 
$k_1 + m_1 = 0$, which leads to $k_0 = k_1$ and $m_0 = m_1$. 
Now we can argue by induction. 
\qed

\Corollary A.2. {\rm([OP95, Lemma~3.1])} If 
the group $\Gamma$ is torsion free, then each 
automorphism of $A$ is graded, i.e., 
$\Aut(A) = \Aut_{\rm gr}(A).$ {\rm(cf.\ Def.~IV.1)} 
\qed

\def\entries{

\[ABFP05 Allison, B. N., Berman, S., Faulkner, J. R., and A. Pianzola, 
{\it Realization of graded-simple algebras as loop algebras}, submitted 

\[AABGP97 Allison, B. N., Azam, S., Berman, S., Gao, Y., and
A. Pianzola, ``Extended Affine Lie Algebras and Their Root Systems,'' 
Memoirs of the Amer. Math. Soc. {\bf 603}, Providence R.I., 1997 

\[BGK96 Berman, S., Gao, Y., and Y. S. Krylyuk, {\it Quantum tori and the 
structure of elliptic quasi-simple Lie algebras}, J. Funct. Anal. {\bf 135} (1996), 
339--389 

\[BL04 Bonahon, F., and L. Xiaobo, {\it Representations of the quantum 
Teichm\"uller space and invariants of surface diffeomorphisms}, 
arXiv:math.GT/0407086v3, 16.11.2004 

\[Bro82 Brown, K. S., ``Cohomology of Groups,'' 
Grad. Texts Math. {\bf 87}, Springer-Verlag, 1987 

\[Br93 Brown, W. C., ``Matrices over Commutative Rings,'' Pure and Appl. Math. 
{\bf 169}, Marcel Dekker, 1993 

\[dCP93 de Concini, C., and C. Procesi, {\it Quantum groups}, in 
``${\cal D}$-modules, representation theory, and quantum groups'', 
Venice (1992), Lecture Notes Math. {\bf 1565} (1993), 31--140

\[Fu70 Fuchs, L., ``Infinite Abelian Groups, Vol. I,'' Pure and Applied Math. {\bf 36}, 
Acad. Press, 1970 

\[GVF01 Gracia-Bondia, J. M., J. C. Vasilly, and H. Figueroa,
``Elements of Non-com\-mu\-ta\-tive Geometry,'' Birkh\"auser Advanced Texts,
Birkh\"auser Verlag, Basel, 2001 

\[Ha00 de la Harpe, P., ``Topics in Geometric Group Theory,'' Chicago Lectures 
in Math., The Univ. of Chicago Press, 2000 

\[Is03 Ismagilov, R. S., {\it The integral Heisenberg group as an infinite 
amalgam of commutative groups}, Math. Notes {\bf 74:5} (2003), 630--636 

\[Jac56 Jacobson, N., ``Structure of Rings,''
Amer. Math. Soc. Coll. Publications {\bf 37}, 1956 

\[KPS94 Kirkman, E., C. Procesi and L. Small, 
{\it A $q$-analog of the Virasoro algebra}, Comm. Alg. {\bf 22:10} (1994), 
3755--3774 

\[La93 Lang, S., ``Algebra,'' 3rd edn., Addison Wesley Publ. Comp., London, 1993

\[New72 Newman, M., ``Integral Matrices,'' Pure and Applied Math. {\bf 45}, 
Acad. Press, New York, 1972 

\[OP95 Osborn, J. M., and D. S. Passman, {\it Derivations of skew polynomial rings}, J. Algebra 
{\bf 176} (1995), 417--448 

\[Pa96 Panov, A. N., {\it Skew fields of twisted rational functions and the skew field 
of rational functions on $\GL_q(n,\K)$}, St. Petersburg Math. J. 
{\bf 7} (1996), 129--143 

}

\def\address
{Karl-Hermann Neeb

Technische Universit\"at Darmstadt 

Schlossgartenstrasse 7

D-64289 Darmstadt 

Deutschland

neeb@mathematik.tu-darmstadt.de}

\references
\lastpage 

\bye